\documentclass[ejs,preprint]{imsart}

\usepackage{amsthm,amsmath,amssymb,mathrsfs,bbm}
\RequirePackage[colorlinks,citecolor=blue,urlcolor=blue]{hyperref}

\RequirePackage{natbib}
\RequirePackage{hypernat}

\doi{10.1214/154957804100000000}
\pubyear{0000}
\volume{0}
\firstpage{0}
\lastpage{0}

\startlocaldefs
\DeclareMathOperator*{\argmin}{arg\,min}
\DeclareMathOperator*{\argmax}{arg\,max}
\renewcommand{\cal}{\mathcal}
\newcommand{\scr}{\mathscr}
\newcommand{\Var}{\mathbb{V}{\rm ar}}
\newcommand{\Cov}{\mathbb{C}{ \rm ov}}
\newcommand{\coloneqq}{\mathrel{\mathop:}=}
\newcommand{\eqqcolon}{=\mathrel{\mathop:}}
\newcommand{\Exp}{\mathbb{E}}
\newcommand{\CExp}[2]{\mathbb{E}\left. \left[ #1 \right| #2 \right] }
\newcommand{\inftynorm}[1]{\left| #1 \right|_\infty}
\newcommand{\norm}[1]{\left\lVert #1 \right\rVert}
\newcommand{\map}[2]{\,{:}\,#1\!\longrightarrow\!#2}
\newcommand{\convP}{ \stackrel{p}{\longrightarrow} }
\newcommand{\convD}{ \stackrel{{\scr D}}{\longrightarrow} }
\newtheorem{definition}{Definition}
\newtheorem{theorem}{Theorem}
\newtheorem{lemma}{Lemma}
\theoremstyle{remark}
\newtheorem{Remark}{Remark}

\endlocaldefs

\begin{document}

\begin{frontmatter}

\title{Kink estimation in stochastic regression with dependent errors and predictors\thanksref{T1}}
\runtitle{Kink estimation in stochastic regression}

\begin{aug}
\author{\fnms{Justin} \snm{Wishart}\corref{}\ead[label=e1]{justin.wishart@sydney.edu.au}}
\address{School of Mathematics and Statistics, F07\\ University of Sydney NSW 2006\\ Australia\\ \printead{e1}}
\and
\author{\fnms{Rafa{\l}} \snm{Kulik}\thanksref{t2}\ead[label=e2]{rkulik@uottawa.ca}}
\address{University of Ottawa\\ \printead{e2}}

\thankstext{T1}{This work was done partially during both authors stay at Banff International Research Station in Banff, Alberta.}
\thankstext{t2}{Research supported by a grant from the Natural Sciences and Engineering Research Coun-
cil of Canada.}

\runauthor{J. Wishart and R. Kulik}
\end{aug}
\begin{abstract}
In this article we study the estimation of the location of jump points in the first derivative (referred to as kinks) of a regression function $\mu$ in two random design models with different long-range dependent (LRD) structures. The method is based on the zero-crossing technique and makes use of high-order kernels. The rate of convergence of the estimator is contingent on the level of dependence and the smoothness of the regression function $\mu$. In one of the models, the convergence rate is the same as the minimax rate for kink estimation in the fixed design scenario with i.i.d. errors which suggests that the method is optimal in the minimax sense.
\end{abstract}




\tableofcontents

\end{frontmatter}

\section{Introduction}

Assume that we observe a bivariate dataset $\left\{ X_i, Y_i \right\}_{i = 1}^n$ that follows the regression model,
\begin{equation}
  Y_i = \mu(X_i) + \sigma(X_i) \varepsilon_i,
  \label{eq:randomregression}
\end{equation}
where $\mu$ is the regression function and $\sigma$ is a deterministic scale function. Also, $\varepsilon_i$ and $X_i$ are the error and random design variables respectively (both being possibly long-range dependent) and $X_i$ has cumulative distribution function $F = F_X \map{\mathbb{R}}{[0,1]}$ that is strictly increasing.

We are interested in testing the presence of a change point in the slope of a regression function $\mu$ and if one exists, estimating its location. We describe this jump in the first derivative of $\mu$ as a kink and denote the change point by $\theta$. Knowledge of this change point will allow us to identify change in trends in the underlying regression function of a non-parametric model. This could explain the change in qualitative or quantitative behaviour of an underlying process.

\subsection{Existing Results}

Before examining the kink estimation under the random design regression model \eqref{eq:randomregression}, we first look at other non-parametric and parametric models and their link to the existing theory for kink point estimation. A change point estimation technique was pioneered by \cite{Goldenshluger-et-al-2006} for estimating change points in the regression function itself, not the kink scenario. The underlying model assumed for their framework was the indirect model with fixed design. The indirect model assumes that the regression function is not observed in practice but a so called `blurred' version of the regression function is observed whereby the regression function has been transformed by a convolution operator. More specifically, the indirect model assumes that observations are realisations of the asymptotic model,
\begin{equation}
	dY(x) = K \, \mu(x) \, dx + \epsilon dB(x).
	\label{eq:asympindirect}
\end{equation}
In the above model the function $K\, \mu (x) = \int_\mathbb{R} K(t-x)\mu(x)\, dx$ represents the convolution of $\mu$ and $K$ and the noise is driven by a regular Brownian motion, $B(x)$ and controlled by $\epsilon \asymp n^{- \frac{1}{2}}$ where the statement $a_n\asymp b_n$ means that the ratio $a_n/b_n$ is bounded above and below by positive constants. The fixed design implies that the design variables $x_i = \frac{i}{n}$ are equally spaced points on the unit interval. The asymptotic model, \eqref{eq:asympindirect} is considered is due to a result by \cite{Brown-Low-1996} that shows \eqref{eq:asympindirect} is asymptotically equivalent to the model,
\begin{equation}
	Y_i = K \, \mu(x_i)  + z_i,
	\label{eq:finiteindirect}
\end{equation}
where $z_i$ is an i.i.d. sequence of error variables. 

The specific estimation technique that \cite{Goldenshluger-et-al-2006} formulated was the zero-crossing technique and it used a particular class kernel functions to identify the change point. Their technique will be adapted for use in this article and is pursued in further detail in \autoref{ssec:approx3rd}. At this stage it will suffice to say that the main result of their paper established that the zero-crossing technique is optimal in the minimax sense under the framework given in \eqref{eq:asympindirect}.

The zero-crossing technique has been applied by \cite{Cheng-Raimondo-2008} to estimate a kink instead of a jump point and was done in the direct model in the fixed design setting. In this framework the observations are assumed to follow a fixed design and realisations derived from the following asymptotic model,
\begin{equation}
	dY(x) = \mu(x) dx + \epsilon dB(x).
 	\label{eq:asympregulardesign}
\end{equation}
Model \eqref{eq:finiteindirect} and their asymptotic equivalents are usually appropriate in practice when a variable is observed at regular intervals indexed by time and the errors are i.i.d. homoscedastic random variables.

More recently, \cite{Wishart-2009} extended the technique further to include long-range dependent (LRD) noise observations instead of independent noise. The kink estimation technique was extended to include the model,
\begin{equation}
	dY(x) = \mu(x) dx + \epsilon^\alpha dB_H(x),
	\label{eq:asympregulardesignLRD}
\end{equation}
where $B_H(x)$ is a fractional Brownian motion with self-similarity index $H \in [\frac{1}{2} , 1)$. The noise process was normalised by $\epsilon^\alpha$ where $\alpha = 2 - 2H$. \cite{Wang-1996} has shown that Model \eqref{eq:asympregulardesignLRD} is the asymptotic equivalent to the discrete model,
\begin{equation}
	y_i = \mu(x_i) + e_i,
	\label{eq:finiteregulardesignLRD}
\end{equation}
where $e_i$ is a LRD sequence of random variables.

In this paper we are interested in model \eqref{eq:randomregression}, which extends the fixed design cases given in models \eqref{eq:finiteindirect}, \eqref{eq:finiteregulardesignLRD} above. They are extended in the sense that the design points are no longer restricted to a uniform grid of points and the scale function  $\sigma(\cdot)$ allows heteroscedasticity for the error terms in the regression model. The analysis of this random design model needs to be considered quite carefully, since the asymptotic behaviour of the estimators will depend on the behaviour of the scale function and on the level of dependence present in the design variables and errors themselves. It has been shown by \cite{Reiss-2008} that there exists an asymptotic equivalence between model \eqref{eq:randomregression} and \eqref{eq:asympregulardesign} when $\sigma(\cdot) \equiv $ constant, and the design variables are independent uniform random variables. However, this is not the case in general. As noted in \cite{Kulik-Raimondo-2009a}, with LRD design variables, model \eqref{eq:randomregression} cannot be equivalent to any asymptotic model, which is in contrast to model \eqref{eq:asympregulardesignLRD} being the asymptotic equivalent to model \eqref{eq:finiteregulardesignLRD} in the fixed design case.

There is an extensive treatment in the literature on both parametric and non-parametric methods for regression models with a random design framework that assume i.i.d. design and error variables. The methodologies used include, but are not limited to, kernel smoothing, wavelet decompositions and orthogonal series. The methods of change point estimation for the random design case have been considered in \cite{Korostelëv-Tsybakov-1993,Gijbels-et-al-1999,Huh-Park-2004}

There is also literature on the fixed design scenario in the presence of long-range dependent errors and the introduction of dependence in the errors always has a detrimental effect on estimation in this scenario. In the context of function estimation some recent treatments of this topic include \cite{Csörgo-Mielniczuk-1995,Wang-1996,Johnstone-Silverman-1997,Johnstone-1999,Cavalier-2004,Kulik-Raimondo-2009a}. For change point estimation work has been done by \cite{Wang-1999,Wishart-2009}.

Then there is a new emerging literature that attempts to combine the two scenarios with random design regression models where the design variables and/or the error variables are LRD. When the framework includes a random design and possibly LRD variables then there is a more subtle asymptotic theory that is based on a delicate balance between the behaviour of the $\sigma$ function and the level of dependence present. This is evident in a current number of papers in the area and will be the case here as well. The interested reader is referred to work by \cite{Robinson-Hidalgo-1997,Guo-Koul-2008} for a parametric linear model approach in this context and to \cite{Csörgo-Mielniczuk-1999,Yang-2001,Mielniczuk-Wu-2004,Kulik-Raimondo-2009} for regression estimation in a non-parametric framework. Finally some studies to estimate change points in the non-parametric context include \cite{Wang-2008,Lin-et-al-2008}.

\subsection{Article Outline}

Some preliminary framework is outlined in \autoref{sec:prelim}, setting up the class of functions that are considered and specific dependence assumptions made in the random design model. The main result of the paper is described in \autoref{sec:mainresult}, along with a brief discussion. The estimation method is explained in detail in \autoref{sec:method}, with a brief outline of the zero-crossing technique in the fixed design and its extension to the random design case. All the necessary proofs of the results are given in \autoref{sec:appendix}.

\section{Preliminaries}
\label{sec:prelim}
\subsection{Smoothness Assumptions and Kernels}
\label{ssec:smoothness}

First we look at the smoothness of the regression function $\mu$ and the properties of the kernel function that was constructed to use the zero-crossing technique by \cite{Cheng-Raimondo-2008}. First we define a class of functions that have domain ${\cal X}$, a kink at $\theta$ and $s \ge 3$ derivatives that exist in the neighbourhood of $\theta$.
\begin{definition}
\label{def:Fs}
We say that $\mu\in {\scr F}_s({\cal X},\theta)$ if, 
\begin{enumerate}
  \item $\mu \map{{\cal X}}{\mathbb{R}}$
  \item $\mu$ has a kink, that is, there exists a  $\theta \in {\cal X}$ and $a_\mu \in \mathbb{R}$ with $a_\mu \neq0$ such that,
  \begin{equation*}
    [ \mu^{(1)}](\theta) = \mu^{(1)}(\theta_+) - \mu^{(1)}(\theta_-) = a_\mu,
  \end{equation*}
  where $\mu^{(1)}(\theta_+)$ and $\mu^{(1)}(\theta_-)$ are the right and left first derivatives of $\mu$ respectively.
  \item The higher order derivatives $\mu^{(i)}$ exist and are finite everywhere and satisfy,
  \begin{equation}
    \mu^{(i)}(\theta_+) = \mu^{(i)}(\theta_-) \qquad \text{for i = 2,3,\ldots,s-1}.
    \label{eq:muderivs}
  \end{equation}
  \item For all $x_+ \in (0,\sup {\cal X}-\theta)$ and $x_- \in(\inf {\cal X} - \theta,0)$,
  \begin{equation}
    \mu^{(1)}(\theta_\pm + x_\pm)- \mu^{(1)}(\theta_\pm) = \sum_{i = 1}^{s-2}\frac{x_\pm^i \mu^{(i+1)}(\theta_\pm)}{i!} + \mathcal{O}(x_\pm^{s-1}). \label{eq:taylor1}
  \end{equation}
  \end{enumerate}
\end{definition}

Condition \textit{4.} should be interpreted in the sense that $\mu^{(1)}$ has a separate Taylor expansion for points to the left and right of $\theta$ respectively. Condition \textit{3.} of \autoref{def:Fs} might seem overly restrictive but is required to exploit the class of Kernel functions that are introduced later in this Section. We will also denote ${\scr F}_s(\theta) = {\scr F}_s(\mathbb{R},\theta)$. For completeness and comparison purposes we will also introduce another smoothness class ${\scr G}_s$ to denote the class of functions that do not have a kink. This class is identical to ${\scr F}_s(\theta)$ except condition 2 and 3 are relaxed in \autoref{def:Fs} in the sense that there does not exist a $\theta \in \mathbb{R}$ such that, $[ \mu^{(1)}](\theta) \neq 0$.

In the fixed design setting, we can assume that the domain of the regression function is $[0,1]$ since any finite interval, $[a,b]$ can be mapped to the $[0,1]$ interval by an affine transformation. However this assumption is not always valid in the general random design case. In particular, if the design variables are LRD then it is required that they have a domain across the whole real line.

To use the zero-crossing technique for this class of regression functions \cite{Cheng-Raimondo-2008} constructed a class of kernel functions via Legendre polynomials and we will denote this class of functions by ${\scr K}_s$. The full description of the zero-crossing technique and the consequent technical details required of the kernel functions are not covered here and the reader is referred to \cite{Goldenshluger-et-al-2006} and \cite{Cheng-Raimondo-2008} respectively for full treatment. However, some key aspects will be given and for our case we will say $K \in {\scr K}_s$, where $s = 2k+1$ and $k \in \mathbb{Z}^+$ if,
\[
	K(x)=K(k,x) =a_k \sum_{j=k-1}^{2k+2}b_{j,k}x^{2j-2k+2}\mathbbm{1}_{[-1,1]}(x)\; ,
\]
where the polynomial coefficients are defined by
\[
a_k\coloneqq \frac{(4k+5)!}{2^{4k+5}(2k)!(2k+2)!}\, ,\quad b_{j,k}\coloneqq \frac{(-1)^{k+j+1}(2j)!}{j!(2k-j+2)!(2j-2k+2)!} .
\]
This class of kernel functions is indexed by the level of smoothness $s$ and is constructed to exploit the extra smoothness of the class ${\scr F}_s(\theta)$. To save on notation we denote $K_i = K^{(i)}$, to represent the $i^{\mathrm th}$ order derivative of $K$.
The kernels have the following properties:
\begin{equation}
  K_i(-1) = K_i(1) = 0 \quad \text{for } i = 1,2,3. \quad \text{and} \quad		K_1(0) = 0.
  \label{eq:Kboundary}
\end{equation}
\begin{align}
\label{eq:K3moment}
  \int_{-1}^1u^j K_3(u)\, du&=0\; , \quad j=0,1,\ldots, 2k.
\end{align}
Property \eqref{eq:K3moment} of ${\scr K}_s$ ensures that the smoothness of ${\scr F}_s\left( \theta\right)$ can be exploited to obtain faster rates of convergence of the estimator $\widehat{\theta}$ in estimating $\theta$. For our purposes of estimation assume that $\mu \in {\scr F}_s(\theta)$ and $\sigma \in {\scr G}_r$ where $s \wedge r \ge 3.$

\subsection{Dependence Assumptions}
\label{ssec:dependence}

Throughout the paper there will be a dependence assumption either among the design random variables or in the error random variables. In particular, the assumed dependence structure is a causal LRD linear process that is defined below.
\begin{definition}
  \label{def:linearprocess}
Let $c_i$ be a set of square summable constant coefficients that are defined,
\begin{equation*}
  c_i \coloneqq \left\{ \begin{array}{rr}
  1, & \text{ if }i = 0,\\
  i^{- (1 + \alpha)/2}L(i), & \text{ if } i \ge 1,
  \end{array}
  \right.\nonumber
\end{equation*}
where $L \map{\mathbb{R}^+}{\mathbb{R}^+}$ is a slowly varying function and $0 < \alpha \le 1$. Then, a random variable $\xi_i$, is said to be a causal LRD linear process if,
\begin{equation*}
  \xi_i = \mu_\xi + \sum_{j = 0}^\infty c_j \eta_{i-j}\nonumber
\end{equation*}
where $|\mu_\xi |< \infty$ and $\eta_i$ are i.i.d. random variables with density $f_\eta$ and moments $\Exp \eta_t = 0$ and $\Exp \eta_t^2 =\left( \sum_{i = 0}^\infty c_j^2 \right)^{-1} \eqqcolon \sigma_\eta^2$. 
\end{definition}
Furthermore, a random variable $\xi_i$ is said to be a causal LRD Gaussian linear process if $\xi_i$ satisfies \autoref{def:linearprocess} and $\left\{ \ldots, \eta_{i-1},\eta_i \right\}$ are i.i.d ${\cal N}\left(0,\sigma_\eta^2\right)$. The case of $\alpha = 1$ is to be interpreted as a short range dependent case and by the construction the random variable has $\Exp \xi_i = \mu_\xi$ and $\Var \xi_i = 1$.  Moreover, it can be shown that $\xi_i$ is a second-order stationary process and has asymptotic covariance structure $\Cov \left( \xi_0 , \xi_k \right) \sim C_0^2 k^{-\alpha}L^2(k)$ where $C_0^2 = \sigma_\eta^2 \int_0^\infty (x^2 + x)^{- (1 + \alpha)/2}\, dx$. Therefore the process exhibits Long-Range Dependence and a consequence of this asymptotic covariance structure is that,
\begin{align}
    \Var\left(\sum_{i=1}^n \xi_i \right)\sim C_1^2n^{2-\alpha}L^2(n),\qquad
    \Var\left(\sum_{i=1}^n \xi_i^2 \right)\sim \left\{\begin{array}{cc}
  C_2^2n^{2-2\alpha}L^4(n), & \quad \text{if }0 < \alpha  < \frac{1}{2},\medskip\\
   C_3^2 n, & \quad \text{if }\frac{1}{2} < \alpha  < 1,
  \end{array}\right.
    \label{eq:asympsumvars}
\end{align}
where $C_1^2 \coloneqq 2C_0^2/((1 - \alpha)(2-\alpha))$, $C_2^2 \coloneqq 4C_0^4/((1-2\alpha)(2-2\alpha))$ and when $1/2 < \alpha < 1$, the covariances $\Cov \left(\xi_0^2,\xi_i^2 \right)$ are summable and $C_3^2 = 1 + 2\sum_{i = 0}^\infty \Cov \left( \xi_0^2, \xi_i^2 \right)$. Also, when $\alpha = 1/2$, $\Var\left(\sum_{i = 1}^n \xi_i^2 \right)$ is asymptotically proportional to a term of order $n$ times another term involving slowly-varying functions. Now throughout the paper, the design variables and error variables are assumed to follow one of the following dependence conditions:
\begin{itemize}
\item[{\rm (A)}] \phantomsection \label{AssumptionA}
The design variables, $\left\{ X_i \right\}_{i = 1}^n$ are i.i.d. random variables with domain ${\cal X}$ and common density $f$ such that $f(x) > 0$ for all $x \in {\cal X}$ and $\sup_{x \in {\cal X}}|f^{(s \wedge r)}(x)| < \infty$. The error variables $\left\{ \varepsilon_i \right\}_{i = 1}^n$ are a causal LRD process with parameter $\alpha_\varepsilon$. Furthermore, the random variables $\left\{ \varepsilon_i\right\}_{i = 1}^n$ are assumed to be independent of $\left\{X_i\right\}_{i = 1}^n$. Under (A), define the associated set of $\sigma$-fields, \[{ \cal G}_i \coloneqq \sigma(\ldots, \eta_{i-1}, \eta_i;X_1, X_2, \ldots, X_i).\]
\item[{\textrm (B)}] \phantomsection \label{AssumptionB}
The design variables, $\left\{ X_i \right\}$ are a causal LRD linear process with parameter $\alpha_x$ where $f_\eta^{(j)}$ is a Lipschitz continuous function for $j = 0,1, \ldots, s$ with $f_X(x) > 0$ for all $x \in \mathbb{R}$. The error variables $\left\{ \varepsilon_i\right\}_{i = 1}^n $, are centred and i.i.d., with a finite variance, independent of $\left\{ X_i\right\}_{i = 1}^n $. Similarly, define the associated set of $\sigma$-fields, \[{\cal F}_i=\sigma(\ldots, \eta_{i-1}, \eta_i;\varepsilon_1,\varepsilon_2,\ldots, \varepsilon_i).\]
\end{itemize}

In both cases, the support of the design variables will be denoted ${\cal X}$. Let $F = F_X$ be the cumulative distribution function of $X$ which is strictly increasing and denote by $F_n(x) = n^{-1} \sum_{i = 1}^n \mathbbm{1}_{\left\{ X_i \le x \right\} }$ the empirical distribution function of $X$. Also let $Q = F^{-1}$ and $Q_n = F_n^{-1}$ be the quantile and empirical quantile functions respectively. We require that $Q$ is Lipschitz, that is, there exists an $L_Q >0$ such that 
\[ 
	|Q(x) - Q(y)| \leq L_Q|x-y|.
\] 
Finally, we need to impose some mild restrictions on $\sigma$. We assume $\sigma$ is bounded away from 0 and $\infty$ in the sense that, 
\[
0 < \inf_{t \in {\cal X}}\sigma(t) < \sup_{t \in {\cal X}}\sigma(t) < \infty\]
 and that $\sigma \in {\scr G}_r$ where $r \ge 3$. Throughout the article we denote by $C$ a general constant that is assumed to be positive and finite but which possibly changes from line to line.

\section{Main Result}
\label{sec:mainresult}

The main result of the paper is concerned with the construction and analysis of an estimator, $\widehat{\theta}$, of the kink location $\theta$. The analysis of the estimator is given in \autoref{thm:rate} and concerns the rate of convergence of $\widehat{\theta}$ to the true the kink location $\theta$. The estimator, $\widehat{\theta}$, will be constructed in \autoref{sec:method} along with the motivations and analysis.

\begin{theorem}
\label{thm:rate}
  Suppose a bivariate sequence of observations $\left\{ X_i , Y_i \right\}$ that follow model \eqref{eq:randomregression} are observed such that $\mu \in {\scr F}_s(\theta)$ and $\sigma \in {\scr G}_r$ where $s\wedge r \ge 3$. Then an estimator, $\widehat \theta$ of the change point, $\theta$, can be constructed such that,
\begin{equation*}
  | \widehat \theta - \theta | = \left\{ \begin{array}{cc}
  {\cal O}_p (n^{- \frac{s}{2s+1}}) ,&\quad \text{under \hyperref[AssumptionA]{Assumption (A)}},\medskip\\
  {\cal O}_p \left(n^{-\frac{s}{2s+1}} \vee \left(n^{- \frac{\alpha_x}{2}}L(n)\right)\right) ,&\quad \text{under \hyperref[AssumptionB]{Assumption (B)}}.
  \end{array}
\right.
\end{equation*}
\end{theorem}
The proof of this Theorem is given at the end of \autoref{sec:method}. The minimax optimality of this result is not pursued in this paper since the lower bounds on the convergence rate of $\widehat \theta$ for the functional class ${\scr F}_s(\theta)$ are not determined in the framework of random design. However, it is worth making the specific point that the obtained rate of convergence under \hyperref[AssumptionA]{Assumption (A)} is the same as the minimax rates for the fixed design case with i.i.d. errors (see \cite{Cheng-Raimondo-2008}). Consequently, it seems reasonable to conjecture that the rates of our estimator are optimal in the minimax sense. 

\section{Kink estimation method}
\label{sec:method}

In this section, the basis of the zero-crossing technique is studied and a brief overview given. Firstly, the zero-crossing technique pioneered by \cite{Goldenshluger-et-al-2006} and applied by \cite{Cheng-Raimondo-2008,Wishart-2009} will be described briefly in \autoref{ssec:approx3rd} and then an adaptation for the random design case constructed in \hyperref[ssec:adapt3rd]{Sections \ref{ssec:adapt3rd} - \ref{ssec:quantile}}. 
\subsection{Approximation of the third derivative for the fixed design model}
\label{ssec:approx3rd}
In the fixed design setting (cf. model \eqref{eq:finiteregulardesignLRD}) it can be assumed without loss of generality that the regression function $\mu$ has domain $[0,1]$. More specifically, assume that $\mu \in {\scr F}_s([0,1],\lambda)$ and estimate $\mu^{(3)}(t)$ by,
\begin{equation*}
  \kappa_h(t) \coloneqq \kappa_h(t,\mu)=h^{-4}\int_0^1 K_3\left( \frac{x-t}{h}\right)\mu(x) \, dx.
\end{equation*}
where $h = h(n)$ is the bandwidth that depends on $n$. Throughout the article it will be assumed that the bandwidth satsifies, at the very least,  $h + \frac{1}{nh} \to 0$, as $n \to \infty$. This is a standard regularity condition for kernel smoothing techniques and additional conditions on the bandwidth will be stated as needed. Using the functional class ${\scr F}_s(\theta)$ and the properties of the kernel function it can be shown that for $t \in (h,1-h)$,
\begin{align}
\kappa_h(t)&=h^{-2}K_{1}\left(\frac{\lambda-t}{h}\right)[\mu^{(1)}](\lambda)+\mathcal{O}(h^{s-3}) =:L_h(t)	+{\cal O}(h^{s-3}),
	\label{eq:1}
\end{align}
where $L_h(t)$ is the localisation term. Indeed, by exploiting the conditions of $K_3$ we can by express $\kappa_h(t)$ as follows. Change variable of integration to obtain,
\begin{align*}
	\kappa_h(t) &= h^{-4} \int_0^1 K_3\left( \frac{x-t}{h}\right)\mu(x)\, dx  \\
		&= h^{-3} \int_{-1}^1 K_3\left( x\right)\mu(t + hx)\, dx.
\end{align*}
The last equality follows because the domain of $K$ is $[-1,1]$ and the values of $t$ are restricted to $t \in (h,1-h)$. This restriction is used to avoid possible edge bias effects from the two sided kernel function. Using integration by parts and exploiting the boundary condition \eqref{eq:Kboundary},
\begin{align}
	\kappa_h(t) &= -h^{-2} \int_{-1}^1 K_2\left( x\right)\mu^{(1)}(t + hx)\, dx.
	\label{eq:khalt}
\end{align}
Let $D = \left\{ t : |\lambda - t| < h\right\}$ and $\tau = (\lambda- t)/h$. Then $|\tau| < 1$ for all $t \in D$. We now split \eqref{eq:khalt} into two integrals,
\begin{equation*}
  \kappa_h(t) = -h^{-2} \int_{-1}^\tau K_2\left( x\right)\mu^{(1)}(t + hx)\,dx -h^{-2} \int_{\tau}^1 K_2\left( x\right)\mu^{(1)}(t + hx)\, dx.
\end{equation*}
To exploit ${\scr F}_s([0,1],\lambda)$ define,
\begin{align*}
  J_h(t) &\coloneqq -h^{-2} \left(  \int_{-1}^\tau K_2\left( x\right)\left( \mu^{(1)}(t + hx) - \mu^{(1)}(\lambda^-)\right)\,dx \right. \\
  & \qquad \qquad \qquad \qquad \left. + \int_{\tau}^1 K_2\left( x\right)\left( \mu^{(1)}(t + hx) - \mu^{(1)}(\lambda^+)\right)\,dx\right) = {\cal O}(h^{s-3}).
\end{align*}
The order bound follows by using \eqref{eq:muderivs} and \eqref{eq:taylor1} in combination with \eqref{eq:K3moment}. Therefore, this allows us to express $\kappa_h(t)$ in the following way,
\begin{align*}
  \kappa_h(t) &= -h^{-2} \int_{-1}^\tau K_2\left( x\right)\mu^{(1)}(\lambda^-)\,dx -h^{-2} \int_{\tau}^1 K_2\left( x\right)\mu^{(1)}(\lambda^+)\,dx + J_h(t) \nonumber \\
  &= h^{-2} K_1(\tau) [\mu^{(1)}](\lambda) + J_h(t)=  L_h(t) + J_h(t).
\end{align*}
This expansion ensures that $\kappa_h(\cdot) = \mathcal{O}(h^{-2})$ for $s \ge 3$, which is assumed to always hold since the third derivative of $\mu$ needs to exist and be finite if it is to be adequately estimated and the method is to make any sense. More specifically we have the following,
\begin{equation}
  \kappa_h(t) = \left\{ \begin{array}{cl}
  {\cal O}(h^{-2}), & \text{if } \mu \in {\scr F}_s(\theta) \text{ and } t\in D\\
  {\cal O}(h^{s-3}), & \text{if $\mu \in {\scr G}_s$} \text{ or } \mu \in {\scr F}_s(\theta) \text{ and } t\notin D.\\
    \end{array}\right.
  \label{eq:orderkappa}
\end{equation}
As seen in all three of the aforementioned papers that use the zero-crossing technique, the $\delta-$separation rate Lemma given below is the technical result that explains why the above representation is effective.
\begin{lemma}[$\delta$-separation rate]
\label{lem:sep}
Let $K \in {\scr K}_s$ and $\mu \in {\scr F}_s([0,1],\theta)$. In what follows the constant $0 < C_q < 1$ depends only the kernel $K_{1} (\cdot)$. Let $h > 0$, $\delta > 0$ be such that $\delta < C_qh$. Let $A_{\delta,h} = \left\{t : \delta < |t - \theta| < C_qh \right\}$. Then for $\kappa_h(t) = \kappa_h(t,\mu)$:
\begin{itemize}
  \item[(a)] $|\kappa_h(\theta)| \leq C h^{s-3}$,
  \item[(b)] for all $t \in A_{\delta,h}$ and $\delta \geq C h^s$, $|\kappa_h(t)| \geq C \delta h^{-3}$,
  \item[(c)] for all $t \in (0,1)$ such that $|\theta - t| > qh$, $|\kappa_h(t)| \leq Ch^{s-3}$.
\end{itemize}
\end{lemma}
The proof of this Lemma is given in \cite{Cheng-Raimondo-2008}. Their proof requires a minor correction as the extra regularity condition $3.$ is needed in the smoothness class ${\scr F}_s(\theta)$.

The main idea of \autoref{lem:sep} allows us to exploit the expansion given in \eqref{eq:1} and focus in on the location of the kink. The kernel function has specific properties to guarantee that a unique global maximum and minimum occurs within order $h$ of the kink point. Furthermore, the estimator was constructed so that the rate of convergence of kink location estimation is minimax for model \eqref{eq:asympregulardesign}. We will seek to adapt these results to the random design setting.

\subsection{Adapted Random Design Estimator of the third derivative}
\label{ssec:adapt3rd}
Now consider $\mu \in {\scr F}_s \left( {\cal X}, \theta \right)$ in model \eqref{eq:randomregression}. An estimator is constructed to exploit the smoothed third derivative of $\mu$ and the argument built around \autoref{lem:sep} discussed in \autoref{ssec:approx3rd}. The most natural extension would be to use the estimator,
\begin{equation}
	\widetilde{\kappa}_h(t) = \frac{1}{nh^4 \widehat{f}_X(t)} \sum_{i = 1}^n K_3 \left( \frac{X_i - t}{h} \right) Y_i,
	\label{eq:randomkappaflaw}
\end{equation}
where $\widehat{f}_X(t)$ is the estimate for the density of $X_i$ at the point $t$ given by,
\begin{equation*}
	\widehat{f}_X(t) = \frac{1}{nh}\sum_{i = 1}^n K \left( \frac{X_i - t}{h} \right).
\end{equation*}
Unfortunately, from a brief computational investigation, the estimator given in \eqref{eq:randomkappaflaw} appears to suffer from poor numerical performance. Instead of using \eqref{eq:randomkappaflaw}, another estimator is constructed by rescaling the design variables by the distribution function $F$ and $\kappa_h(t)$ is estimated in the random design setting by,
\begin{align}
	\widehat{\kappa}_h(t) &= \frac{1}{nh^4} \sum_{i = 1}^n Y_i K_{3} \left( \frac{F(X_i) - t}{h}\right)
	\label{eq:khat}.
\end{align}
This estimator was chosen since it also is a proxy for the fixed design estimator given in \autoref{ssec:approx3rd} and seems to exhibit better numerical performance than \eqref{eq:randomkappaflaw}. The estimator given in \eqref{eq:khat} is also an unbiased estimate of the smoothed third derivative,
\begin{align}
	\Exp \widehat{\kappa}_h(t) &= h^{-4} \Exp \mu(X_1) K_{3} \left( \frac{F(X_1) - t}{h}\right)\nonumber\\
		&= h^{-4} \int_\mathbb{R} \mu(u) K_{3} \left( \frac{F(u) - t}{h}\right)\, dF(u)\nonumber\\
		&= h^{-4} \int_0^1 \mu_F(x) K_{3} \left( \frac{x - t}{h}\right)\, dx = \kappa_h(t,\mu_F),
	\label{eq:2}
\end{align}
where $\mu_F(\cdot) = \mu(Q(\cdot))$. If $\mu \in {\scr F}_s(\theta)$, then $\mu_F \in {\scr F}_s([0,1],\lambda)$ where $\theta = Q(\lambda)$. In \eqref{eq:2}, the observed quantity is the smoothed third derivative of $\mu_F$, which, coupled with \autoref{lem:sep} and the argument shown in \autoref{ssec:approx3rd} is equivalent to estimating a kink location $\lambda$ for the function $\mu_F$ in the fixed design setting.

Therefore with the above argument, an estimator $\widehat{\theta}$ of a kink location of the regression function $\mu$ in the random design setting is constructed that is approximately the same as the estimator for kink location $\lambda$ of $\mu_F$ in the fixed design setting. This is done by estimating the value of $\lambda$ by $\widehat{\lambda}$ using the established zero-crossing technique in the fixed design setting and then rescaling $\lambda$ back by the quantile function to obtain an estimate of $\theta$. Thus to assess the performance of our estimator we need to check that the convergence of $\widehat{\kappa}_h(t)$ to $\kappa_h(t)$ is sufficiently fast. To do this consider the two following processes,
\begin{equation}
\begin{split}
\begin{aligned}
	\gamma_i(t) &= \mu(X_i) K_3\left( \frac{F(X_i) - t}{h}\right)\\
	\zeta_i(t) &= \sigma(X_i) K_3\left( \frac{F(X_i) - t}{h}\right).\label{eq:gammazeta}
\end{aligned}
\end{split}
\end{equation}
With these definitions, the overall accuracy of the estimator can be decomposed into, 
\begin{equation}
	\widehat{\kappa}_h(t)=  \kappa_h(t) + b_h(t)+ Z_h(t) ,
  \label{eq:kappahatdecomp}
\end{equation}
where $b_h(t)$ and $Z_h(t)$ represent the respective stochastic error and stochastic bias contributions to the estimator and are given by, 
\[
  b_h(t) = n^{-1} h^{-4} \sum_{i=1}^n \left( \gamma_i(t) - \Exp \gamma_1(t) \right), \qquad Z_h(t) = n^{-1} h^{-4}\sum_{i=1}^n \zeta_i(t) \varepsilon_i.
\]
The analysis of the above terms are given in the next subsection. 

\subsection{Probabilistic Behaviour for the Adapted Estimator}
\label{sssec:probbounds}

In this section the analysis of the stochastic bias and stochastic error terms are considered before proceeding to the next stage of the zero-crossing technique to ensure that the stochastic contributions do not overwhelm the signal generated by the $\kappa_h(t)$ term. The proofs of the claims in this section will be deferred to \autoref{sec:appendix}.

The first term to be considered is the stochastic bias term which did not appear in previous kink analyses pursued by \cite{Cheng-Raimondo-2008,Wishart-2009} since there is some stochastic contribution by adapting the fixed design estimator to the the random design framework. Therefore, this term needs to be appropriately dealt with and the next Lemma is a useful tool that considers this term.

\begin{lemma}
\label{LIL}
Consider a function $\mu \map{{\cal X}}{\mathbb{R}}$ such that $\mu'$ exists and is bounded. Then define the function \[ \gamma^*_i(t) = \left( \mu(X_i) - \mu_F(t) \right) K_3\left( \frac{F(X_i) - t}{h} \right) .\] 
If the design variables follow \hyperref[AssumptionA]{Assumption (A)} then,
\begin{equation*}
  \sup_{t \in (0,1)} \left| \sum_{i = 1}^n \left(  \gamma^*_i(t)- \Exp \gamma^*_i(t)\right)\right|= o_{\rm a.s. }\left( \sqrt{nh^3 \left| \log h\right|}\right).
\end{equation*}
If the design variables follow \hyperref[AssumptionB]{Assumption (B)} then,
\begin{equation*}
  \sup_{t \in (0,1)} \left| \sum_{i = 1}^n \left(  \gamma^*_i(t)- \CExp{\gamma^*_i(t) }{{\cal F}_{i-1}}\right)\right| = {\cal O}_p\left( \sqrt{nh^3 \left| \log h\right|}\right).
\end{equation*}
\end{lemma}
Note that the two claims in given in \autoref{LIL} follow from the uniform law of iterated logarithms for independent variables and an similar iterated logarithm result for martingale difference sequences.
 
We now state some central and non-central limit theorems for the estimator, $\widehat{\kappa}_h(t)$. The convergence of the estimator $\widehat{\kappa}_h(t)$ under both \hyperref[AssumptionA]{Assumption (A)} and \hyperref[AssumptionB]{(B)} is contingent on the size of the bandwidth relative to the level of dependence $\alpha$. The specific details of this relationship between $h$ and $n^\alpha$ will be shown in detail inside the Theorems. Roughly speaking, if the bandwidth is too `large' compared to $\alpha$ then the dependence of the random variables dominate and the estimator converges to a process that needs to be normed by a sequence that relies on $\alpha$. Conversely, if the bandwidth is `small' compared to $\alpha$ then the dependence of the random variables is negligible and a regular central limit theorem holds with a norming sequence that is not reliant on $\alpha$. In the forthcoming Theorems the extra smoothness of the regression and variance functions are exploited to be able to obtain an estimator that is not as sensitive to the level of dependence. In practice, this extra level of smoothness will most likely be unknown. Due to its common occurrence in the subsequent Theorems, define the asymptotic variance term, $\upsilon^2(t) \coloneqq \left( \sigma_F^2(t) + \mu_F^2(t)\right) \int_{-1}^1 K_3^2 (x)\, dx$. The following Theorem deals with the case of \hyperref[AssumptionA]{Assumption (A)}.

\begin{theorem} 
\label{CLTA}
Let $K \in {\scr K}_{s \wedge r}$, $\mu \in {\scr F}_s$, $\sigma \in {\scr G}_r$ with $s \wedge r \ge 3$ and $t \in (h,1-h)$. Also if the design variables and error random variables follow \hyperref[AssumptionA]{Assumption (A)} and the bandwidth $h = h(n)$ also satisfies,
\begin{equation}
  h^{2(s \wedge r)+1} n^{1-\alpha_\varepsilon}L^2(n) \rightarrow 0\quad \text{as } n \to \infty,
  \tag{A1}
\label{eq:bandA1}
\end{equation}
 then the following convergence result holds,
\begin{align}
  \sqrt{nh^7} \left( \widehat{\kappa}_h(t) - \kappa_h(t) \right) &\stackrel{{\scr D}}{\longrightarrow} {\cal N}\left( 0, \upsilon^2(t)\right).
   \label{eq:CLTA1}
\end{align}
Conversely, if the bandwidth $h = h(n)$  satisfies,
\begin{equation}
	 h^{2(s \wedge r)+1} n^{1-\alpha_\varepsilon}L^2(n) \rightarrow \infty\quad \text{as } n \to \infty,
	\tag{A2}
  \label{eq:bandA2}
\end{equation}
then,
\begin{align*}
  \frac{n^{\frac{\alpha_\varepsilon}{2}} h^{3-(s \wedge r)}}{L(n)} \left( \widehat{\kappa}_h(t) - \kappa_h(t) \right) &\stackrel{{\scr D}}{\longrightarrow} {\cal N}\left( 0, C_1^2\upsilon_*^2(t)\right).
\end{align*}
where \[ \upsilon_*(t) = \frac{\sigma_F^{(s \wedge r)}(t) }{(s \wedge r)!} \int_{-1}^1  x^{s \wedge r}K_3(x)\, dx .\]
\end{theorem}

\autoref{CLTB1} and \autoref{CLTB2} deal with case under \hyperref[AssumptionB]{Assumption (B)} and give the central limit theorems when there is a `small' or `large' bandwidth respectively. In the `large' bandwidth scenario a stronger assumption is used whereby the design variables are a causal LRD Gaussian linear process.

\begin{theorem} \label{CLTB1} Let $K \in {\scr K}_{s \wedge r}$, $\mu \in {\scr F}_s$, $\sigma \in {\scr G}_r$ with $s \wedge r \ge 3$ and $t \in (h,1-h)$. If the design variables and error random variables follow \hyperref[AssumptionB]{Assumption (B)} and the bandwidth $h = h(n)$ satisfies,
\begin{align}
 h^{7} n^{1-\alpha_x}L^2(n) \to 0 \quad \text{as } n \to \infty,\tag{B1}
 \label{eq:bandB1}
\end{align}
then the estimator obeys the following law,
\begin{equation}
	\sqrt{nh^7} \left( \widehat{\kappa}_h(t) - \kappa_h(t)\right)\stackrel{{\scr D}}{\longrightarrow}  {\cal N}\left(0,\upsilon^2(t)\right).
	\label{eq:regularCLT}
\end{equation}
\end{theorem}
\begin{theorem} 
\label{CLTB2}
Let $K \in {\scr K}_{s \wedge r}$, $\mu \in {\scr F}_s$, $\sigma \in {\scr G}_r$ with $s \wedge r \ge 3$ and $t \in (h,1-h)$. Assume the design variables and error random variables follow \hyperref[AssumptionB]{Assumption (B)} and that the design variables are a causal LRD Gaussion linear process. If the bandwidth $h =h (n)$ satisfies,
\begin{align}
h^{7} n^{1-\alpha_x}L^2(n) \to \infty\quad \text{as }n  \to \infty,\tag{B2}
 \label{eq:bandB2}
\end{align}
and the estimator  $\widehat{\kappa}_h(t)$ has a Hermite rank of 1 then the the estimator obeys the following law,
\[ 
  \frac{n^{\frac{\alpha_x}{2}}}{L(n)} \left( \widehat{\kappa}_h(t) - \kappa_h(t) \right) \convD {\cal N}(0,C_1^2{\cal H}_1(t))
\]
where 
\[ {\cal H}_1(t) = \frac{\kappa_h(t)}{s_X^3\sigma_\eta\phi\left( \Phi^{-1} ( t)\right)}\int_\mathbb{R} \phi\left(\frac{\Phi^{-1}\left( t\right)-u}{s_X}\right) \left(\Phi^{-1}\left( t\right) - u\right)\phi\left( \frac{u}{\sigma_\eta}\right) \, du,\]
 $s_X^2 = 1 - \sigma_\eta^2$ and $\phi$ and $\Phi$ are the standard normal density and cumulative distribution functions respectively.
\end{theorem}

\begin{Remark}
\label{rem:CLTB2}
If the estimator $\widehat{\kappa}_h(t)$ has Hermite rank $q$ for some $q \in \left\{2,3,\ldots \right\}$ then the asymptotic distribution depends on the size of the bandwidth relative to $q\alpha$. Firstly, if $n^{1-q\alpha_x}h^7 L^{2q} (n) \to \infty$ then it can be shown using a similar argument used in the Proof of \autoref{CLTB2} with the result of Theorem 2 of \cite{Avram-Taqqu-1987} that the normed process $n^{q\alpha_x/2}L^{-q}(n)\left(\widehat{\kappa}_h(t) - \kappa_h(t)\right) \convD {\cal H}_q(t){\scr H}_q$ where, 
\[ {\cal H}_q(t) = \frac{\kappa_h(t)}{s_X^2\sigma_\eta\phi\left( \Phi^{-1} ( t)\right)}\int_\mathbb{R} \phi\left(\frac{\Phi^{-1}\left( t\right)-u}{s_X}\right) H_q\left(\frac{\Phi^{-1}\left( t\right) - u}{s_X}\right)\phi\left( \frac{u}{\sigma_\eta}\right) \, d u \] and $H_q(x)$ is the Hermite polynomial of degree $q$ and ${\scr H}_q$ is the Hermite-Rosenblatt process,
\[ {\scr H}_q = \sqrt{\frac{q! (1-q\alpha)(2 - q \alpha)}{2\left( \int_0^\infty (x^2 + x)^{-(1 + \alpha)/2}\, dx\right)^q}} \int_{-\infty < x_1 < x_2 \ldots < x_q < 1} \left\{ \int_0^1 \prod_{i = 1}^q \left( \left( y - x_i \right)^+ \right)^{-\frac{\alpha+1}{2}}dy \right\}dB(x_1) \ldots dB(x_q)\]
where $B$ denotes a standard Brownian motion. In \cite{Avram-Taqqu-1987}, they considered Appell polynomials for a generalised sequence of stationary LRD random variables. In our case the LRD variables are Gaussian and consequently the Appell polynomials reduce to the Hermite polynomials. On the other hand, if the bandwidth satisfies $n^{1-q\alpha_x}h^7 L^{2q} (n) \to 0$ then \eqref{eq:regularCLT} holds.
\end{Remark}

As will be seen in \autoref{ssec:kinkdetect}, some large deviations results are needed to be able to to distinguish between the signal generated by the $\kappa_h(t)$ term and the stochastic bias and noise contributions. Unfortunately, a slightly weaker large deviations result is proved under \hyperref[AssumptionA]{Assumption (A)} in \autoref{EVTA}. In particular we assume that the scale function, $\sigma(\cdot) \equiv \sigma$, is constant however this restriction could possibly be relaxed by using a different method. The large deviations result for \hyperref[AssumptionB]{Assumption (B)} in \autoref{EVTB} does not carry this restriction and the scale function need not be constant. 
\begin{theorem}
\label{EVTA}
Let $K \in {\scr K}_{s\wedge r}$ and the design and error variables satisfy \hyperref[AssumptionA]{Assumption (A)}. Further assume that the bandwidth $h = h(n)$ also satisfies,
\begin{equation}
   \frac{\left| \log h\right|^3}{nh^3} + \frac{L^2(n)\left| \log h\right|^2}{n^{\alpha_\varepsilon}h^{\frac{4}{3}}} \to 0 \quad \text{as }n \to \infty.
  \label{eq:EVTAband}
\end{equation}
Then define,
\begin{align*}
  S_n^A(t) \coloneqq \sum_{i = 1}^n \Psi_i(t) = \frac{1}{\sqrt{nh\, \upsilon^2(t)}}\sum_{i = 1}^n K_3\left( \frac{F(X_i) - t}{h}\right)\left(\mu_F(t) + \sigma\varepsilon_i \right).\nonumber
\end{align*}
Also define, 
\begin{equation}
  B_n(x) = \sqrt{2 \log n } + \frac{x}{\sqrt{2 \log n }} - \frac{1}{\sqrt{2 \log n }} \left( \frac{1}{2} \log \log n + \log \left( 2 \sqrt{\pi} \right)\right)
  \label{eq:Bn}
\end{equation}
and a partition of $[0,1]$,
\begin{equation}
  T_n = \left\{ t_j = 2 h j, j = 1, \ldots, m_n - 1\right\}
  \label{eq:Tn}
\end{equation}
where $m_n = \lceil \frac{1}{2h} \rceil.$ Then,
\begin{equation*}
  \lim_{n \to \infty} P\left( \sup_{t \in T_n} \left|S_n^A(t)\right| \le B_{m_n}(x)\right) = e^{-2e^{-x}},
\end{equation*}
for all $x \in \mathbb{R}$.
\end{theorem}

\begin{theorem}
\label{EVTB}
Let $K \in {\scr K}_{s\wedge r}$ and the design and error variables satisfy \hyperref[AssumptionB]{Assumption (B)} and assume that the bandwidth $h = h(n)$ also satisfies,
\begin{equation}
   \frac{\left| \log h\right|^3}{nh^3} + \frac{L^2(n)\left| \log h\right|^2}{n^{\alpha_x}h^{\frac{4}{3}}} + h^{2\left( s \wedge r\right)+1}n \to  0 \quad \text{as } n \to \infty.
  \label{eq:EVTBband}
\end{equation}
Then define,
\begin{align*}
  S_n^B(t) \coloneqq \sum_{i = 1}^n \Xi_i(t) = \frac{1}{\upsilon(t)\sqrt{nh}}\sum_{i = 1}^n  K_3\left(\frac{F(X_i) - t}{h}\right) \left( \mu_F(t) +  \sigma(X_i) \varepsilon_i \right).\nonumber
\end{align*}
with $B_n(x)$ and $T_n$ defined by \eqref{eq:Bn} and \eqref{eq:Tn} respectively, then,
\begin{equation*}
  \lim_{n \to \infty} P\left( \sup_{t \in T_n} \left|S_n^B(t) \right|\le B_{m_n}(x)\right) = e^{-2e^{-x}},
\end{equation*}
for all $x \in \mathbb{R}$.
\end{theorem}

\subsection{Localisation Step}

Recall from \eqref{eq:1}, that the probe function given by $\kappa_h(t)$ gives a signal from the localisation term, $L_h(t)$ with some approximation error and the estimator adds a stochastic bias and error term,
\begin{equation}
	\widehat\kappa_h(t) = L_h(t)+{\cal O}(h^{s-3}) + Z_h(t) + b_h(t).
	\label{eq:Okappa}
\end{equation}
Clearly, $h^{-2} > h^{s-3}$, since $s \ge 3$. So to be able to discern the signal generated from $L_h(t) = {\cal O}(h^{-2})$, it is required that $L_h(t)$ dominates the stochastic terms, $Z_h(t)$ and $b_h(t)$. 

By construction of the Kernel function, (cf. \cite{Cheng-Raimondo-2008}), $K_{1}(\cdot)$ has a unique minimum and maximum in the interval $[-1,1]$, so that $K_{1}(\cdot/h)$ has a unique minimum and maximum in the interval of a length $\mathcal{O}(h)$. Consequently, $L_h(\cdot)$ has a unique extrema near $t^*=\theta+\mathcal{O}(h)$ and $t_*=\theta-\mathcal{O}(h)$. As in the fixed design scenario considered by \cite{Cheng-Raimondo-2008,Wishart-2009} define,
\[
	t_*\coloneqq \argmin_{t \in (0,1)} L_h(t)\; ,\qquad t^* \coloneqq \argmax_{t \in (0,1)} L_h(t).
\]
However, in practice the location of $t_*$ and $t^*$ are not known and estimated using $\widehat\kappa_h(t)$ with,
\[
\widehat t_*=\argmin_{t \in (0,1)} \widehat \kappa_h(t)\; ,\qquad \widehat t^*=\argmax_{t \in (0,1)} \widehat \kappa_h(t).
\]
If $\mu_F \in {\scr F}_s([0,1],\lambda)$ then, 
\begin{equation}
  |L_h(t^*)+L_h(t_{*})|\ge C h^{-2}.
  \label{eq:Lhextrema}
\end{equation}
There are two respective bandwidth restrictions, (\eqref{eq:bandA1}, \eqref{eq:bandA2}; \eqref{eq:bandB1}, \eqref{eq:bandB2}) for the asymptotic behaviour of the estimator under each of the \hyperref[AssumptionA]{Assumption (A)} and \hyperref[AssumptionB]{Assumption (B)} respectively. Starting with \eqref{eq:bandA1} and \eqref{eq:bandB1}, to have a well defined signal, it is required that, $h^{-2} \geq C n^{-\frac{1}{2}}h^{-\frac{7}{2}} \Rightarrow h \ge C n^{-\frac{1}{3}}$. Furthermore, since it is assumed that $s \wedge r \ge 3$, to ensure that \eqref{eq:CLTA1} and \eqref{eq:regularCLT} always hold it suffices to choose $h$ such that $h \le Cn^{- \frac{1}{7} + (\alpha_x \vee \alpha_\varepsilon) /7 - \delta}$, for some $\delta> 0$ or,
\begin{equation}
	C n^{-\frac{1}{3}+\delta} < h < Cn^{- \frac{1}{7} - \delta}
	\label{eq:optimalh}
\end{equation}
 for some $\delta > 0$. With this choice, the bandwidth restrictions given by \eqref{eq:bandA1} and \eqref{eq:bandB1} will always hold. 

It is worth noting that under this choice, the order of the stochastic terms does not involve $\alpha_x$ or $\alpha_\varepsilon$, the level of dependence. Note that $h$ is chosen in a very similar manner if $\varepsilon_i$ and $X_i$, $i\ge 1$, were i.i.d. Consequently, there will be no influence of the (long range) dependence on the change point estimation. The influence of the long range dependence will only affect testing purposes of the threshold used to determine if a signal is genuine and this will be discussed in the next subsection.

\subsection{Kink Detection step}
\label{ssec:kinkdetect}
For simplicity in notation, assume that $[\mu_F]^{(1)}(\lambda) > 0$, which means,  $t_* < t^*$ (a similar argument follows if $[\mu_F]^{(1)}(\lambda) < 0 \Rightarrow t_* > t^*$.) To detect a kink, first standardise the statistic $\widehat{\kappa}_h(t)$ to have unit variance. This will allow us to appropriately notice if there is a change-point present when the observed extrema of $\widehat{\kappa}_h(t)$ exceed the threshold for the noise process. Define this standardised process as,

\begin{equation}
	\mathcal{T}_{\widehat{\kappa}}(t) \coloneqq \frac{\sqrt{nh^7 }\widehat{\kappa}_h(t)}{\upsilon(t)} . 
	\label{eq:Tstat}
\end{equation}
Then by \eqref{eq:Okappa} and \eqref{eq:Tstat} the $\mathcal{T}_{\widehat{\kappa}}(t)$ process has expansion,
\begin{align}
		\mathcal{T}_{\widehat{\kappa}}(t) &=\frac{n^{\frac{1}{2}} h^{\frac{7}{2}}}{\upsilon(t)} L_h(t) + o(n^\frac{1}{2}h^{\frac{1}{2}}) + \frac{n^{\frac{1}{2}} h^{\frac{7}{2}}}{\upsilon(t)}\left( Z_h(t) + b_h(t) \right).
	\label{eq:Tstatexp}
\end{align}

As seen earlier, the information regarding a kink is generated by the $L_h(t)$ process. A thresholding regime will be considered to be able to distinguish between the signal generated by $L_h(t)$ against the noise signal generated by the $Z_h(t)$ and $b_h(t)$ terms. This thresholding will be split into the two scenarios for \hyperref[AssumptionA]{Assumption (A)} and \hyperref[AssumptionB]{(B)}. 

Begin by firstly giving a general decomposition of the estimator for both cases by using, $\gamma^*_i(t) = \left( \mu(X_i) - \mu_F(t) \right)K_3\left(\frac{F(X_i) - t}{h}\right) = \gamma_i(t) + \mu_F(t) K_3 \left( \frac{F(X_i) - t}{h}\right)$ and using \eqref{eq:gammazeta} and \eqref{eq:kappahatdecomp}. So,
\begin{align} 
  \mathcal{T}_{\widehat{\kappa}}(t) &= \frac{\sqrt{nh^7}}{\upsilon(t)}\kappa_h(t) + \frac{1}{\upsilon(t)\sqrt{nh}} \sum_{i = 1}^n \left( \gamma_i(t) - \Exp \gamma_1(t)  + \zeta_i(t) \varepsilon_i \right)\nonumber\\
  &= \frac{\sqrt{nh^7}}{\upsilon(t)}\kappa_h(t) + \frac{1}{\upsilon(t)\sqrt{nh}} \sum_{i = 1}^n \left( \gamma^*_i(t) - \Exp \gamma^*_1(t) + \mu_F(t) K_3 \left( \frac{F(X_i) - t}{h}\right)  + \zeta_i(t) \varepsilon_i \right).
  \label{eq:EVTdecomp}
\end{align}

First assume $\sigma(X_i) \equiv \sigma$, constant, and focus on \hyperref[AssumptionA]{Assumption (A)}. By an application of \autoref{LIL} and \eqref{eq:K3moment},
\begin{align*}
  \mathcal{T}_{\widehat{\kappa}}(t) &= \frac{\sqrt{nh^7}}{\upsilon(t)}\kappa_h(t) + o_{\rm a.s. }\left( \left|\log h\right| \right) + S_n^A(t).
\end{align*}
From \autoref{EVTA}, it is known that $S_n^A(t)$ will diverge to infinity no faster than $\sqrt{2 \left| \log 2h \right|}$. Also, if $\mu \in {\scr G}_s$, then from \eqref{eq:orderkappa}, $\kappa_h(t) = {\cal O}(h^{s-3})$ and
\begin{equation}
\lim_{n \to \infty} P \left( \sup_{t \in T_n} \mathcal{T}_{\widehat{\kappa}}(t) \ge \sqrt{2 \left| \log 2h \right|} \right)=0.
\label{eq:Pthreshold}
\end{equation}
However, if $\mu \in {\scr F}_s\left(\theta\right)$, then \eqref{eq:Lhextrema} holds and by \eqref{eq:Tstatexp}, $\max_{t \in \left( t_*, t^*\right)}{\cal T}_{\widehat \kappa}(t) \ge C n^{\frac{1}{2}}h^{\frac{3}{2}} > \sqrt{2 \left| \log 2h \right|}$ and a kink is detected when, 
\begin{equation}
  \max_{t \in T_n} \left| {\cal T}_{\widehat \kappa}(t) \right| \ge \sqrt{2 \left| \log 2h \right|}.
\label{eq:kinkdetect}
\end{equation}

A very similar argument holds for \hyperref[AssumptionB]{Assumption (B)}. In this case assume that the scale function $\sigma \in {\scr G}_r$ with $r \ge 3$ and proceed as before. In conjunction with \eqref{eq:EVTdecomp} and \eqref{eq:K3moment} apply \autoref{LIL},
\begin{align*}
  \mathcal{T}_{\widehat{\kappa}}(t) &= \frac{\sqrt{nh^7}}{\upsilon(t)}\kappa_h(t)  + S_n^B(t) + {\cal O}_p\left( h\sqrt{\left|\log h\right|} \right)\nonumber\\
  &= \frac{\sqrt{nh^7}}{\upsilon(t)}\kappa_h(t)+ S_n^B(t) + o_p\left( 1 \right).
\end{align*}
The bandwidth restriction \eqref{eq:optimalh} guarantees that \eqref{eq:EVTBband} and consequently \autoref{EVTB} holds. Then for \hyperref[AssumptionB]{Assumption (B)} the same argument applies that was used to show \eqref{eq:Pthreshold} for \hyperref[AssumptionA]{Assumption (A)}.

This thresholding technique does raise some restrictions that could possibly be removed by another technique. Recall from \eqref{eq:optimalh}, that $h > Cn^{-\frac{1}{3}+\delta}$ for some $\delta > 0$ is required to be able to distinguish the signal from the stochastic terms. Also, \eqref{eq:EVTAband} and \eqref{eq:EVTBband} are required to be able to apply \autoref{EVTA} and \autoref{EVTB} respectively and obtain a large deviation result for the process. Therefore to ensure both conditions are satisfied, it is sufficient to consider $\alpha_x > \frac{4}{9}$ or $\alpha_\varepsilon > \frac{4}{9}$.

\subsection{Zero Crossing Technique}

The idea behind the zero-crossing technique is that within the interval $\widehat{A}_h = [\widehat t_*, \widehat t^*]$, $\widehat{\kappa}_h(t) \approx \kappa_h(t)$. Using \autoref{lem:sep} we can locate the zero-crossing-time of $\kappa_h(t)$ which occurs at $t=\lambda$ with an accuracy of order $\delta, \, \delta < h$. First minimise $|\widehat{\kappa}_h(t)|$ within the interval $\widehat{A}_h$:
\begin{equation*}
	\widehat{\lambda} = \argmin_{t \in \widehat{A}_h} |\widehat{\kappa}_h(t)| = \argmin_{t \in \widehat{A}_h} |\mathcal{T}_{\widehat{\kappa}}(t)|.
\end{equation*} 
	By comparing \eqref{eq:1} with the bounds in \autoref{lem:sep} we see that the minimum is well defined if,
	\begin{equation}
		\delta h^{-3} \geq C h^{s-3} \quad \text{and} \quad \delta h^{-3} \geq C n^{-\frac{1}{2}}h^{-\frac{7}{2}}.
	\label{eq:Tbound}
\end{equation}
	We will obtain the best possible accuracy if we choose $\delta$ as small as possible, as long as both inequalities of \eqref{eq:Tbound} still hold. The left hand expression of \eqref{eq:Tbound} implies that $\delta \asymp h^s$ and substituting this into the right hand expression of \eqref{eq:Tbound} we derive the order of the smallest possible bandwidth 
\begin{equation*}
		h_* \asymp n^{-\frac{1}{2s + 1}}.
\end{equation*}
	 We now apply \autoref{lem:sep} with $\delta_* = h_*^s$ to locate the change point $\lambda$ in $\mu_F$ with an accuracy of order,
\begin{equation*}
		\left| {\widehat \lambda} - \lambda \right|= \delta_{*} = h_*^s \asymp n^{-\frac{s}{2s+1}}.
\end{equation*}

\subsection{Modified Estimator of Kink}
\label{ssec:quantile}

Recall that $\theta = Q(\lambda)$. In practice the true distribution function $F$ is unknown, so it is estimated in the usual manner by the empirical distribution function $F_n(x)= n^{-1}\sum_{i = 1}^n \mathbbm{1}_{\left\{ X_i \le x\right\}}$ and consequently can obtain an estimator of $Q$ via the empirical quantile function $Q_n(\cdot).$ Estimate $\theta$ by, $\widehat \theta = Q_n(\widehat \lambda)$. The rate of convergence of this estimator is evaluated below,
\begin{align}
	|\widehat{\theta} - \theta| &= |Q_n(\widehat \lambda) - Q(\lambda)|\nonumber\\
		&\leq |Q_n(\widehat \lambda) - Q(\widehat \lambda)| + |Q(\widehat \lambda) - Q(\lambda)|\nonumber\\
		&\leq |Q_n(\widehat \lambda) - Q(\widehat \lambda)| + L_Q |\widehat \lambda - \lambda| \nonumber\\
		&\le |Q_n(\widehat \lambda) - Q(\widehat \lambda)| + \mathcal{O}_p(n^{- \frac{s}{2s+1}}).
	\label{eq:genrate}
\end{align}
The rate of convergence in \eqref{eq:genrate} is therefore contingent on the maximum of the rate from the generalised quantile process for the design variables or the rate from the initial unscaled kink estimator. Under \hyperref[AssumptionA]{Assumption (A)}, the quantile process involves independent and identically distributed design variables and for all $t \in (0,1)$, 
\begin{equation}
 \left|Q_n(t) - Q(t) \right| = {\cal O}_p(n^{-\frac{1}{2}})
\label{eq:QA}
\end{equation}
(see \cite{Csörgo-1983} and references therein for a detailed treatment). For \hyperref[AssumptionB]{Assumption (B)}, the rate is dependent on $\alpha_x$ and for all $t \in (0,1)$,
\begin{equation}
 \left|Q_n(t) - Q(t) \right| = {\cal O}_p(n^{-\frac{\alpha_x}{2}}L(n))
\label{eq:QB}
\end{equation}
(see Theorem 5.1 of \cite{Ho-Hsing-1996}). Therefore, using \eqref{eq:QA} and \eqref{eq:QB} in \eqref{eq:genrate},
\begin{equation*}
  | \widehat \theta - \theta | = \left\{ \begin{array}{cc}
  {\cal O}_p (n^{- \frac{s}{2s+1}}) ,&\quad \text{under \hyperref[AssumptionA]{Assumption (A)}}.\medskip\\
  {\cal O}_p \left(n^{-\frac{s}{2s+1}} \vee \left(n^{- \frac{\alpha_x}{2}}L(n)\right)\right) ,&\quad \text{under \hyperref[AssumptionB]{Assumption (B)}},
  \end{array}
\right.
\end{equation*}
which proves \autoref{thm:rate}.

\begin{Remark}
The method can be extended to the multiple kink scenario by observing multiple instances of \eqref{eq:kinkdetect}. For each instance of \eqref{eq:kinkdetect} there is a corresponding interval $\widehat{A}_h$ and the localisation and zero-crossing-time steps are executed on each of those intervals to produce an estimate for each kink location. The interested reader is referred to \cite{Cheng-Raimondo-2008,Wishart-2009} for a more detailed treatment of the method in the multiple kink scenario.

\end{Remark}

\section{Mathematical Appendix}
\label{sec:appendix}

Before giving the proofs, some notation is described. Let $X$ denote a random variable and denote the $L_p$-norm $\norm{X}_p^p = \Exp \left|X\right|^p$ and $\norm{\cdot} = \norm{\cdot}_2$. For a function $f \map{{\cal X}}{\mathbb{R}}$ denote the sup-norm $\inftynorm{f} = \sup_{x \in {\cal X}} \left| f(x)\right|$. Throughout this Section a Taylor expansion of composite functions will be used to exploit the vanishing moment condition of $K_3$. For the Taylor expansion to be well defined, the derivatives of the composite functions need to exist. A generalised chain rule for composite functions exists (see the Fa\`a di Bruno formula from \cite{Encinas-et-al-2005} and references therein), and these are of the form,
\begin{equation}
 \frac{d^n}{dx^n} f(g(x)) \coloneqq \frac{d^n}{dx^n} (f \circ g)(x)=   \sum_{\mathbf{k} \in {\cal K}_n} \frac{n!}{k_1! k_2! \ldots k_n!} (f^{(k)} \circ g)(x) \prod_{i = 1}^n \left(\frac{g^{(i)}(x)}{i!}\right)^{k_i}
\label{eq:FaadiBruno}
\end{equation}
where ${\cal K}_n = \left\{ k_i \in \left\{\mathbb{Z}^+ \cup {0}\right\}: k_1 + 2k_2 + \ldots + nk_n = n\right\}$ and $k = \sum_{i = 1}^n k_i$. Also, through tedious but elementary calculus it can be shown that, the $n^{\rm th}$ derivative of $Q = F^{-1}$ will exist, and the Taylor expansions of $ \mu_F$ and $\sigma_F$ up to order $n$ will exist if $f^{(n)}$ exists.

\begin{proof}[Proof of \autoref{LIL}]
  Begin with the proof of the first claim under \hyperref[AssumptionA]{Assumption (A)}. Since $\varsigma_i(t)$ will be non-zero only if $F(X_i) \in (t-h,t+h)$, there exists a $\tau_i \in (-1,1)$ that depends on $X_i$ such that,
\begin{align*}	
  \gamma^*_i(t) &=K_3 \left( \frac{F(X_i) - t}{h}\right)\left( \,\mu_F(t + \tau_i h) - \mu_F(t)\right)\nonumber\\
  &= h \tau_i K_3 \left( \frac{F(X_i) - t}{h}\right)\mu_F^{(1)}(t + h \xi_i \left|\tau_i\right|h) \eqqcolon h \nu_i(t),
\end{align*}
and $\xi_i$ depends on $\tau_i$. The $\nu_i(t)$ terms are independent random variables, each of which have variance that is of order $h$. Therefore by the Law of Iterated Logarithm (see \cite{Bingham-1986}) we have the following result,
\begin{equation*}
  \limsup_{n \to \infty} \frac{1}{\sqrt{nh \log \log n}} \sum_{i = 1}^n \Big( \nu_i(t) - \Exp \nu_i(t)\Big) = -\liminf_{n \to \infty} \frac{1}{\sqrt{nh \log \log n}} \sum_{i = 1}^n \Big( \nu_i(t) - \Exp \nu_i(t)\Big)= C \qquad \text{a.s.,}
\end{equation*}
Therefore we have,
\begin{align*}
  \sum_{i = 1}^n \Big( \gamma^*_i(t) - \Exp \gamma^*_i(t)\Big) &= h \sum_{i = 1}^n \Big( \nu_i(t) - \Exp \nu_i(t)\Big)\nonumber\\
  &= {\cal O}_{\text{a.s.}}\left(\sqrt{nh^3\log \log n}\right)\\
  &= o_{\text{a.s.}}\left(\sqrt{nh^3\left|\log h\right|}\right)
\end{align*}
which proves the first claim of the Lemma. Now to concentrate on the claim for \hyperref[AssumptionB]{Assumption (B)}, a proof of a similar claim in Lemma 4 of \cite{Zhao-Wu-2006} is adapted to our framework. This technique bounds the martingale difference sequence $\gamma^*_i(t) - \CExp{\gamma^*_i(t)}{{\cal F}_{i-1}}$ above and below by two discretised martingale difference sequences and uses an exponential martingale inequality to gain the required probabilistic bounds. To do this, again exploit the Taylor expansion of $\mu$ in \autoref{def:Fs} and use the fact that  $Support(K_3) = [-1,1]$, which means that there exists a $\tau_i$ dependent on $X_i$ with $\left| \tau_i \right| \le 1$ such that $F(X_i) = t + \tau_i h$ and,
\begin{align}
  \gamma^*_i(t) &= K_3 \left( \frac{F(X_i) - t}{h}\right) \left( \,\mu_F(t + \tau_i h) - \mu_F(t) \right) \mathbbm{1}_{\left( t-h,t+h\right)}\left(F(X_i)\right)\nonumber\\
  &= \tau_i h K_3 \left( \frac{F(X_i) - t}{h}\right)  \mu_F^{(1)}(t + \xi \left| \tau_i \right| h)  \mathbbm{1}_{\left(t-h,t+h\right)}\left(F(X_i)\right),
\label{eq:gammamuidecomp}
\end{align}
where $\left| \xi \right| \le 1$. Then split the function in \eqref{eq:gammamuidecomp} into its positive and negative parts by defining $\xi_i \coloneqq t + \xi\left| \tau _i \right| h$ and
$\tau_i \mu_F^{(1)}(\xi_i) = \left(\tau_i \mu_F^{(1)}(\xi_i)\right)^+- \left(\tau_i \mu_F^{(1)}(\xi_i)\right)^-\eqqcolon \mu_{F,1}^+(\xi_i) - \mu_{F,1}^-(\xi_i)$ where $f^+ = (f\vee0)$, $f^- = (-f\wedge 0)$ denote the respective positive and negative parts of $f$. Then,
\begin{align}
  \gamma^*_i(t) &= h \left( \mu_{F,1}^+(\xi_i) - \mu_{F,1}^-(\xi_i)\right)\left( K_3^+\left( \frac{F(X_i)-t}{h}\right) - K_3^-\left( \frac{F(X_i)-t}{h}\right)\right) \mathbbm{1}_{(t-h,t+h)}((F(X_i)))\nonumber\\
  &=  h \left( \mu_{F,1}^+(\xi_i)K_3^+\left( \frac{F(X_i)-t}{h}\right) - \mu_{F,1}^+(\xi_i)K_3^-\left( \frac{F(X_i)-t}	{h}\right)\right.\nonumber\\
  &\left.\qquad \qquad  - \mu_{F,1}^-(\xi_i)K_3^+\left( \frac{F(X_i)-t}{h}\right) + \mu_{F,1}^-(\xi_i)K_3^-\left( \frac{F(X_i)-t}{h}\right)\right) \mathbbm{1}_{(t-h,t+h)}((F(X_i)))\nonumber\\
  &\coloneqq \left( \varsigma^{++}_i(t)-\varsigma^{+-}_i(t)-\varsigma^{-+}_i(t)+\varsigma^{--}_i(t) \right),
\label{eq:gammadecomp}
\end{align}
By the linearity of the conditional expectation operator and \eqref{eq:gammadecomp} we can decompose the martingale difference sequence into parts,
\begin{align}
  \gamma^*_i&(t) - \CExp{\gamma^*_i(t)}{{\cal F}_{i-1}}\nonumber\\
  &=  \varsigma^{++}_i(t)- \CExp{\varsigma^{++}_i(t)}{{\cal F}_{i-1}}-\left(\varsigma^{+-}_i(t)-\CExp{\varsigma^{+-}_i(t)}{{\cal F}_{i-1}}\right)\nonumber\\
  &\qquad \qquad - \left(\varsigma^{-+}_i(t)-\CExp{\varsigma^{-+}_i(t)}{{\cal F}_{i-1}}\right)+  \varsigma^{--}_i(t)-\CExp{\varsigma^{--}_i(t)}{{\cal F}_{i-1}}
  \label{eq:gammaMDSdecomp}
\end{align}
To begin with we will concentrate on the first martingale difference term on the RHS of \eqref{eq:gammaMDSdecomp} and bound it above and below by a discretised version that does not depend on $t$ directly. For this discretization let $N = \lceil \left(nh^{-3}\right)^{\frac{1}{2}} \rceil$ and $t_j = \frac{j}{N}$ where $0 \le j \le N.$ Then for any $t \in [0,1]$ there exists a $j$ such that $t \in [t_j, t_{j+1})$ and the distance $|t_{j+1} - t_j| = {\cal O}(N^{-1})$. Define the two new tweaked martingale difference sequences versions of $\varsigma_i^{++}(t)$,
\begin{align*}
  \overline{\varsigma}_{i,j}^{++} &= h\mu_{F,1}^+(\xi_i) \frac{1}{2}\left\{ K_3^+ \left( \frac{F(X_i) - t_j}{h}\right) + K_3^+ \left( \frac{F(X_i) - t_{j+1}}{h}\right) \right\}\mathbbm{1}_{(t_{j}-h,t_{j+1}+h)}(F(X_i))\\
  \underline{\varsigma}_{i,j}^{++} &= h\mu_{F,1}^+(\xi_i) K_3^+ \left( \frac{F(X_i) - t_j}{h}\right)\mathbbm{1}_{(t_{j+1}-h,t_j+h)}(F(X_i))
\end{align*}
It can be shown that, the martingale difference sequence $\varsigma_i^{++}(t) - \CExp{\varsigma_i^{++}(t)}{{\cal F}_{i-1}}$ can be bounded uniformly in $t$ above and below by,
\begin{align*}
	\underline{\varsigma}_{i,j}^{++}  - \CExp{\underline{\varsigma}_{i,j}^{++}}{{\cal F}_{i-1}}-CN^{-1} \le \varsigma_i^{++}(t) - \CExp{\varsigma_i^{++}(t)}{{\cal F}_{i-1}} &\le \overline{\varsigma}_{i,j}^{++}  - \CExp{\overline{\varsigma}_{i,j}^{++}}{{\cal F}_{i-1}}+CN^{-1}.
\end{align*}
We have the following result,
\begin{align*}
  \sup_{t \in (0,1)} \left| \sum_{i= 1}^n \left(\varsigma_i^{++}(t) -\CExp{\varsigma_i^{++}(t)}{{\cal F}_{i-1}} \right)\right| &\le \max_{0 \le j \le N-1} \left( \left| \overline{S_n}(j) \right| +\left| \underline{S_n}(j) \right| \right) + CnN^{-1}\nonumber\\
  &= \max_{0 \le j \le N-1} \left( \left| \overline{S_n}(j) \right| +\left| \underline{S_n}(j) \right| \right) + o\left( \sqrt{nh^3 \left| \log h\right|}\right),
\end{align*}
where for each fixed $j$, $\underline{S_n}(j)$ and $\overline{S_n}(j)$ are martingales with respect to the filtration ${\cal F}_n$ and are defined,
\begin{align*}
  \underline{S_n}(j) = \sum_{i = 1}^n \left( \underline{\varsigma}_{i,j}^{++} - \CExp{\underline{\varsigma}_{i,j}^{++}}{{\cal F}_{i-1}} \right)\\
  \overline{S_n}(j) = \sum_{i = 1}^n \left( \overline{\varsigma}_{i,j}^{++} - \CExp{\overline{\varsigma}_{i,j}^{++}}{{\cal F}_{i-1}} \right).
\end{align*}
These martingales will be bounded by an exponential martingale inequality. Consider firstly the martingale $\underline{S_n}(j)$, its martingale differences are bounded $\left|\underline{\varsigma}_{i,j}^{++} - \CExp{\underline{\varsigma}_{i,j}^{++}}{{\cal F}_{i-1}}\right|\le 2h\inftynorm{\mu^{(1)}}\inftynorm{K_3} \eqqcolon C_b h $. Also using the Lipschitz property of $Q$ and the bounded domain of $K_3$,
\begin{align*}
  \CExp{\left( \underline{\varsigma}_{i,j}^{++} - \CExp{\underline{\varsigma}_{i,j}^{++}}{{\cal F}_{i-1}}\right)^2}{{\cal F}_{i-1}} &\le \int_\mathbb{R} \left( h \mu_{F,1}^+(u) K_3^+\left( \frac{F(u)-t_j}{h}\right) \right)^2\left. f_X( u \right| {\cal F}_{i-1}) \, du \\
  &\le 2 h^3  L_Q \inftynorm{K_3}^2 \inftynorm{\mu_{F,1}^+} \inftynorm{f_\eta} \eqqcolon C_{cv}h^3.
\end{align*}
Then, a martingale inequality for bounded differences given by Theorem 1.5A of \cite{de-la-Pena-1999} can be used to yield,
\begin{align}
	P \left(  \underline{S_n}(j) \ge x \right) &\le \exp\left\{ -\frac{x}{2a}  \sinh^{-1} \left(\frac{ax}{2y}\right) \right\},
	\label{eq:delaPena}
\end{align}
where $a = C_bh$ and $y = C_{cv} nh^3$. Furthermore if $ax/2y = o(1)$ then using a Taylor expansion of $\sinh^{-1}$,
\begin{equation}
  \sinh^{-1} \left( \frac{ax}{2y}\right) = \frac{ax}{2y} + o\left( \left(\frac{ax}{2y}\right)^2\right).
  \label{eq:Taylorsinh}
\end{equation}
Now consider the chance that $\max_{1 \le j \le n} \underline{S_n}(j) $ exceeds the threshold $x = C_T\sqrt{nh^3\left| \log h \right|}$ for some $C_T>0$ which combined with $a = C_bh$ and $y = C_{cv}nh^3$ implies, $ax/2y = {\cal O}\left(\sqrt{\left| \log h\right|/nh}\right) = o(1)$ and by \eqref{eq:delaPena} and \eqref{eq:Taylorsinh},
\begin{align}
  P \left( \underline{S_n}(j) \ge C_T \sqrt{nh^3 \left| \log h \right|} \right) &\le \exp \left\{ - \frac{C_T^2}{4C_{cv}} \left| \log h \right| + o\left( 1\right)  \right\}
  \label{eq:underSPbound}
\end{align}
So, fix $\epsilon > 0$ and use \eqref{eq:underSPbound},
\begin{align}
  P \left( \max_{0 \le j \le N-1 } \underline{S_n}(j)  \ge C_T \sqrt{nh^3\left| \log h\right|}\right) &\le P \left( \bigcup_{j = 0}^{N-1} \left\{   \underline{S_n}(j) \ge C_T\sqrt{nh^3\left| \log h\right|} \right\}\right) \nonumber\\
  &\le \sum_{j = 0}^{N-1} P \left( \underline{S_n}(j) \ge C_T\sqrt{nh^3\left| \log h\right|} \right) \nonumber\\
  &\le N \exp \left\{ - \frac{C_T^2}{4C_{cv}} \left| \log h \right|\right\}\exp \left\{ o\left( 1\right)  \right\}\nonumber\\
  &\le Cn^{\frac{1}{2}}h^{C_T^2/4C_{cv}- \frac{3}{2}}.
  \label{eq:PboundubarGj1}
\end{align}
 By choosing $C_T$ large enough will ensure that $Cn^{\frac{1}{2}}h^{C_T^2/4C_{cv}- \frac{3}{2}}< \epsilon$. The similar conclusion can be reached that for any $\epsilon > 0$ there exists a finite constant $C$ such that,
\begin{equation}
  P \left( -\max_{0 \le j \le N-1 } \underline{S_n}(j)  \ge C \sqrt{nh^3\left| \log h\right|}\right)  < \epsilon.
  \label{eq:PboundubarGj2}
\end{equation}
Therefore, \eqref{eq:PboundubarGj1} and \eqref{eq:PboundubarGj2} ensure that,
\begin{equation*}
  \max_{0 \le j \le N-1 } \underline{S_n}(j)={\cal O}_p(\sqrt{nh^3\left| \log h\right|}).
\end{equation*}
Using a comparable argument, the same conclusion can be reached for the $\overline{S_n}(j)$,
\begin{equation*}
  \max_{0 \le j \le N-1 } \left| \overline{S_n}(j)\right|={\cal O}_p(\sqrt{nh^3\left| \log h\right|}).
\end{equation*}
Also, a similar technique can be used to bound the other martingale difference terms given in \eqref{eq:gammaMDSdecomp}, details omitted.
\end{proof}

\begin{proof}[Proof of \autoref{CLTA}]
To prove the Theorem we appeal to similar results that were shown by \cite{Wu-Mielniczuk-2002,Kulik-2008} by decomposing the stochastic terms into two parts, a martingale part and a LRD part. This is done by defining, 
\begin{align*}
	\chi_i(t) = \frac{\left(\zeta_i(t) - \Exp \zeta_1(t)\right) \varepsilon_i + \gamma_i(t) - \Exp \gamma_1(t)}{\sqrt{n \, \left( \Var \zeta_1(t) + \Var \gamma_1(t) \right)}}
\end{align*}
and then decomposing the standardised estimator $\widehat{\kappa}_h(t)$ into two terms,
\begin{align}
	\sqrt{nh^7} \left( \widehat{\kappa}_h(t) - \kappa_h(t) \right) &= \sqrt{nh^7} \left( Z_h(t) + b_h(t) \right)\nonumber\\
		&= \frac{1}{\sqrt{nh}} \left( \sum_{i=1}^n \zeta_i(t) \varepsilon_i  + \sum_{i= 1}^n \left( \gamma_i(t) - \Exp \gamma_1(t) \right) \right)\nonumber\\
		&= \sqrt{h^{-1}\, \left( \Var \zeta_1(t) + \Var \gamma_1(t) \right)} \sum_{i = 1}^n \chi_i(t) +\frac{\Exp \zeta_1(t)}{\sqrt{nh}} \sum_{i = 1}^n \varepsilon_i.\label{eq:CLTXidecomp}
\end{align}
The Theorem will follow by showing that either the first or last term on the RHS of \eqref{eq:CLTXidecomp} dominates under the bandwidth conditions \eqref{eq:bandA1} or \eqref{eq:bandA2} respectively. More specifically, it will be shown that the dominating term will follow a CLT and the other term converges to zero in probability; then Slutsky's Theorem completes the proof. Firstly consider the case where \eqref{eq:bandA1} holds, then apply the martingale CLT of \cite{Brown-1971} to show,
\begin{equation}
	\sum_{i=1}^n \chi_i(t) \stackrel{{\scr D}}{\longrightarrow} {\cal N} (0 , 1).
	\label{eq:CLTXi}
\end{equation}
Note that $\left\{ \chi_i(t), {\cal G}_i \right\}$ form a martingale difference sequence. So it remains to check that the sum of the conditional variances converge in probability to the unconditional sum and the Lindeberg condition holds. Before we prove the Lindeberg condition note that for $t \in (h,1-h)$,
\begin{align}
	\Exp \zeta_1^2(t) 	&= \int_\mathbb{R} \sigma^2(x) K_3^2\left(\frac{F(x)-t}{h}\right) \, dF(x)  = h\int_{-1}^1 \sigma_F^2(t + hu) K_3^2\left(u\right) \, du.
	\label{eq:orderzeta21}
\end{align}
Exploiting \eqref{eq:K3moment} and the assumption that $\sigma \in {\scr G}_r$,
\begin{align}
	\Exp \zeta_i(t) &= h\int_{-1}^1 \sigma_F(t + hu) K_3\left(u\right) \, du\nonumber\\
		&= \frac{h^{(s \wedge r)+1}}{(s \wedge r)!} \int_{-1}^1 \sigma_F^{(s \wedge r)}(t + \tau hu) u^{s \wedge r} K_3\left(u\right) \, du = h^{(s \wedge r)+1}\upsilon_*(t), 
	\label{eq:orderexpzeta}
\end{align}
where $\tau \in (0,1)$. Therefore, using \eqref{eq:orderzeta21} and \eqref{eq:orderexpzeta},
\begin{align*}
	\Var \zeta_1(t) &= h \int_{-1}^1 \sigma_F^2(t + hu) K_3^2\left(u\right) \, du - \frac{h^{2(s \wedge r)+2}}{\left((s \wedge r)!\right)^2}\left( \int_{-1}^1 \sigma_F^{(s \wedge r)}(t + \tau hu) u^{s \wedge r} K_3\left(u\right) \, du \right)^2.
\end{align*}
Due to the fact that the bandwidth is assumed to follow $h \in (0,1)$, there exists a $h_0$  such that for all $ 0 < h \le h_0$,
\begin{equation}
	 \Var \zeta_1(t) \ge \frac{h \inf_{x \in \mathbb{R}}\left|\sigma^2(x)\right| }{2} \int_{-1}^1 K_3^2\left(u\right) \, du.
	\label{eq:orderzeta2}
\end{equation}
From \eqref{eq:orderexpzeta}, it follows, $h^{-\frac{1}{2}} \Exp \zeta_1(t) = o(1)$ and from \eqref{eq:orderzeta21}, $h^{-1} \Exp \zeta_1^2(t)  \to  \sigma_F^2(t) \int_{-1}^1 K_3^2 \left( u\right) \, du$. Therefore, $h^{-1} \Var \zeta_1(t) = h^{-1} \left( \Exp \zeta_1^2(t) - \left( \Exp \zeta_1(t) \right)^2\right)\to  \sigma_F^2(t) \int_{-1}^1 K_3^2 \left( u\right) \, du.$
Also, the same argument applies for the $\gamma_i(t)$ term to yield,
\begin{equation}
  h^{-1} \left(  \Var \zeta_1(t) + \Var \gamma_1(t) \right) \stackrel{h \to 0}{\longrightarrow} \upsilon^2(t).
  \label{eq:varkappaCLT}
\end{equation}
Now the Lindeberg condition is shown to hold. Let $\epsilon > 0$ be arbitrary,
\begin{align}
  \sum_{i=1}^n \Exp \chi_i^2 (t) \mathbbm{1}_{\left\{ |\Xi_i(t)| > \epsilon \right\}} &= n\Exp \chi_1^2 (t) \mathbbm{1}_{\left\{ |\chi_1(t)| > \epsilon \right\}} \nonumber\\
  &=  \frac{\Exp \left[ \left( \varepsilon_1\left(\zeta_1(t) - \Exp \zeta_1(t) \right) + \gamma_1(t) - \Exp \gamma_1(t) \right)^2 \mathbbm{1}_{A_n}\right]}{ \Var \zeta_1(t) + \Var \gamma_1(t)}.
  \label{eq:lindebergcondA1}
\end{align}
where $A_n = \left\{ \left| \varepsilon_1 \left(\zeta_1(t) - \Exp \zeta_1(t)\right) + \gamma_1(t) - \Exp \gamma_1(t) \right| > \epsilon\sqrt{n \, \left( \Var \zeta_1(t) + \Var \gamma_1(t) \right)} \right\}$. The size of this set can be maximised using \eqref{eq:orderzeta2}, 
\begin{align}
  A_n & \subset \left\{ 2 \inftynorm{K_3} | \varepsilon_1| \left(  \inftynorm{\sigma} + \inftynorm{\mu}\right) > \epsilon\sqrt{n  \Var \zeta_1(t) } \right\} \nonumber\\
  & \subset \left\{ 2 \inftynorm{K_3}| \varepsilon_1| \left(  \inftynorm{\sigma} + \inftynorm{\mu}\right) > \epsilon\sqrt{n h \frac{ \inf_{x \in \mathbb{R}}\left|\sigma^2(x)\right| }{2} \int_{-1}^1  K_3^2\left(u\right) \, du } \right\}.
  \label{eq:lindebergcondA2}
\end{align}
Using the fact that $nh \to \infty$ and $h \to 0$ as $n \to \infty$ we see that $A_n \to \emptyset$, the empty set. Consequently with \eqref{eq:lindebergcondA1}, \eqref{eq:lindebergcondA2} and $n\Exp \chi_1^2(t) < \infty$ imply that,
\begin{align*}
  \sum_{i=1}^n \Exp \chi_i^2 (t) \mathbbm{1}_{\left\{ |\chi_i(t)| > \epsilon \right\}} \stackrel{n \to \infty}{\longrightarrow} 0,
\end{align*} 
and the Lindeberg condition holds. By a consequence of \eqref{eq:asympsumvars}, let $\epsilon > 0$ be arbitrary,
\begin{align*}
P\left( \left|\frac{1}{n}\sum_{i = 1}^n \varepsilon_i \right| > \epsilon \right) &\le \frac{1}{n^2 \epsilon^2} \Var\left( \sum_{i = 1}^n \varepsilon_i \right) \le \frac{C_1^2 n^{-\alpha} L^2(n)}{\epsilon^2}= o(1),\\
P\left( \left|\frac{1}{n}\sum_{i= 1}^n \varepsilon_i^2  - 1 \right| > \epsilon\right) &\le \frac{1}{n^2 \epsilon^2} \Var \left( \sum_{i = 1}^n \varepsilon_i^2 \right) \le \frac{\left( C_2^2n^{-1} \vee C_3^2n^{-2\alpha}L^2(n) \right)}{\epsilon^2} = o(1).
\end{align*}
Then by the above, the sum of the conditional variances to converge in probability to one:
\begin{align*}
\sum_{i = 1}^n \CExp{\chi_i^2(t)}{{\cal G}_{i-1}} &= \frac{\sum_{i = 1}^n \varepsilon_i^2 \,\Var \zeta_1(t) + n \, \Var \gamma_1(t) + 2 \Cov \left( \zeta_1(t), \gamma_1(t) \right)\sum_{i = 1}^n \varepsilon_i }{n \, \left( \Var \zeta_1(t) + \Var \gamma_1(t) \right)} \convP 1,
\end{align*}
and by the martingale CLT, \eqref{eq:CLTXi} follows.

Now we show that the last term on the RHS of \eqref{eq:CLTXidecomp} converges in probability to zero. Consider an arbitrary $\epsilon > 0$, then using \eqref{eq:orderexpzeta} and \eqref{eq:asympsumvars}, 
\begin{align}
	P \left( \left| \frac{\Exp \zeta_1(t)}{\sqrt{nh}}  \sum_{i = 1}^n \varepsilon_i \right|  > \epsilon \right) &\le \frac{ \left( \Exp \zeta_1(t)\right)^2 }{\epsilon^2 n h}\Var \left( \sum_{i = 1}^n \varepsilon_i \right)\nonumber\\
		&\le Ch^{2(s \wedge r)+1} n^{1 - \alpha} L^2(n)\nonumber\\
  &= o(1),
\label{eq:convPXiAend}
\end{align}
and the last line follows by the bandwidth restriction given in \eqref{eq:bandA1}. Thus, the proof of the first claim under the `small' bandwidth scenario holds.

Consider now the `large' bandwidth scenario. Using \eqref{eq:CLTXidecomp}, \eqref{eq:CLTXi} and \eqref{eq:orderexpzeta}, 
\begin{align}
 {\widehat \kappa_h(t)} - \kappa_h(t) &= {\cal O}_p\left( n^{-\frac{1}{2}}h^{-\frac{7}{2}}\right) + \frac{\upsilon_*(t)}{nh^{3-(s \wedge r)}}\, \sum_{i = 1}^n \varepsilon_i.
  \label{eq:CLTA2decomp}
\end{align}
Also, from \cite{Ho-Hsing-1997}, it is known that
\begin{equation}
  \frac{1}{n^{1 - \frac{\alpha}{2}}L(n)}\sum_{i = 1}^n \varepsilon_i \convD {\cal N}(0,C_1^2).
  \label{eq:HoHsing}
  \end{equation}
Therefore, normalising the expression on \eqref{eq:CLTA2decomp},
\begin{align*}
  \frac{n^{\frac{\alpha}{2}}h^{3- (s \wedge r)}}{L(n)} \left({\widehat \kappa_h(t)} - \kappa_h(t)\right) &= {\cal O}_p\left( h^{- \frac{1}{2} - (s \wedge r)} n^{- \frac{1 - \alpha}{2}}L^{-1}(n) \right) + \frac{\upsilon_*(t)}{n^{1- \frac{\alpha}{2}}L(n)} \, \sum_{i = 1}^n \varepsilon_i,
\end{align*}
and the result follows from \eqref{eq:bandA2} and \eqref{eq:HoHsing} with Slutsky's Theorem
\end{proof}

\begin{proof}[Proof of \autoref{CLTB1}]

First break down the estimator into its separate martingale and LRD part in a similar fashion to the method employed in the proof of \autoref{CLTA}. 
Using \eqref{eq:EVTdecomp}, apply \autoref{LIL},
\begin{align}
\widehat{\kappa}_h(t) - \kappa_h(t)&= \frac{1}{nh^4} \sum_{i = 1}^n \Big( \gamma_i(t) - \Exp \gamma_i(t)  + \zeta_i(t) \varepsilon_i\Big)\nonumber\\
  &= {\cal O}_p\left( \sqrt{\frac{\left| \log h\right|}{nh^5}} \right) + \frac{1}{nh^4} \sum_{i = 1}^n \Big(\CExp{ \gamma^*_i(t)}{{\cal F}_{i - 1}}- \Exp \gamma^*_i(t)\Big)\nonumber\\
  &\qquad \qquad  + \frac{1}{nh^4} \sum_{i = 1}^n \Big(\mu_F(t)K_3 \left( \frac{F(X_i) - t}{h}\right) + \zeta_i(t) \varepsilon_i\Big)\nonumber\\
  &= \frac{1}{nh^4} \sum_{i = 1}^n \left(\mu_F(t)\left(K_3 \left( \frac{F(X_i) - t}{h}\right) - \CExp{K_3\left( \frac{F(X_i) - t}{h}\right)}{{\cal F}_{i-1}}\right) + \zeta_i(t) \varepsilon_i\right)\nonumber\\
  &\qquad \qquad  + \frac{1}{nh^4} \sum_{i = 1}^n \Big(\CExp{ \gamma_i(t)}{{\cal F}_{i - 1}}- \Exp \gamma_i(t) \Big) + {\cal O}_p\left( \sqrt{\frac{\left| \log h\right|}{nh^5}} \right)
  \label{eq:CLTBdecomp}
\end{align}
Define the standardised stochastic terms, 
\[
  \Delta_i(t) \coloneqq \frac{\zeta_i(t)\varepsilon_i + \mu_F(t) \left( K_3 \left( \frac{F(X_i) - t}{h} \right) - \CExp{K_3 \left( \frac{F(X_i) - t}{h} \right)}{{\cal F}_{i-1}}\right)}{\upsilon(t)\sqrt{nh }}.
\]
Then in a similar fashion to the Proof of \autoref{CLTA} it will be shown by the martingale CLT of \cite{Brown-1971} that,
\begin{equation}
	\sum_{i = 1}^n \Delta_i(t) \convD {\cal N}(0,1).
	\label{eq:CLTB}
\end{equation}
Indeed, $\Delta_i(t)$ is a martingale difference sequence with respect to the $\sigma$-fields $\left\{ {\cal F}_i \right\}$. Thus we need to check that the Lindeberg condition holds and that the sum of the conditional variances converge in probability to 1. First, focus on the convergence of the conditional variances. The conditional variances can be broken into two parts,
\begin{align}
  \sum_{i = 1}^n \CExp{\Delta_i^2(t)}{{\cal F}_{i-1}} &= \frac{\mu_F^2(t)}{nh\,\upsilon^2(t)}\sum_{i = 1}^n  \CExp{\left(K_3\left( \frac{F(X_i)-t}{h}\right)-\CExp{K_3\left( \frac{F(X_i)-t}{h}\right)}{{\cal F_{i-1}}}\right)^2}{{\cal F_{i-1}}}\nonumber\\
  	&\qquad \qquad \qquad +\sum_{i = 1}^n\frac{ \CExp{\zeta_i^2(t)}{{\cal F}_{i-1}}}{nh\,\upsilon^2(t)}.
  \label{eq:condvarB}
\end{align}
Dealing with the second term on the RHS of \eqref{eq:condvarB}, use Lemma 1 of \cite{Zhao-Wu-2008},
\begin{align}
  \frac{1}{nh}\sum_{i = 1}^n\CExp{\zeta_i^2(t)}{{\cal F}_{i-1}} &= \frac{1}{nh}\sum_{i = 1}^n \Exp \zeta_i^2(t) + \frac{1}{nh}\sum_{i = 1}^n\Big( \CExp{\zeta_i^2(t)}{{\cal F}_{i-1}} - \Exp \zeta_i^2(t) \Big) \nonumber\\
  &= \int_{- \frac{t}{h}}^{\frac{1-t}{h}} \sigma_F^2(t+hx)K_3^2(x) \, dx+{\cal O}_p(n^{- \frac{\alpha}{2}}L(n))\nonumber\\
  &= \sigma_F^2(t)\int_{- 1}^1 K_3^2(x) \, dx + {\cal O}(h^2)+ {\cal O}_p(n^{- \frac{\alpha}{2}}L(n))
  \label{eq:condvarB1}
\end{align}
To bound the first term of \eqref{eq:condvarB}, a bound is required for $\CExp{K_3 \left( \frac{F(X_i)-t}{h}\right)}{{\cal F}_{i-1}}^2$. Define $X_{i,i-1} \coloneqq X_i - \eta_i = \mu_X + \sum_{j = 1}^\infty c_j \eta_{i-j}$ and $Z_i \coloneqq s_X^{-1}(X_{i,i-1} - \mu_X)$ and define $\widetilde{f}_\eta(x) \coloneqq f_X\left( x \big| {\cal F}_{i-1}\right) = f_\eta (x - X_{i,i-1})$ and $g(x) = 1/x$. Then $X_{i,i-1}$ and $Z_i$ are ${\cal F}_{i-1}$-measurable and for all $t \in (h,1-h)$ the conditional expectation can be evaluated as follows.
\begin{align}
   \Exp \bigg[ K_3\left( \frac{F(X_i)-t}{h}\right)\bigg|{\cal F}_{i-1}\bigg] &= \int_\mathbb{R}  K_3\left( \frac{F(v)-t}{h}\right)f_X\left( v \big| {\cal F}_{i-1}\right)\, dv\nonumber\\
  &= h\int_{-1}^1 K_3\left( x\right)\left( \widetilde{f}_\eta \circ Q\right)(t + hx)\left(g \circ f_X \circ Q\right)(t+hx)\, dx.
\label{eq:CExpK31}
\end{align}
Use a Taylor expansion of the composite functions, $p(t) \coloneqq \left(\widetilde{f}_\eta \circ Q \right)(t)$ and $q(t) \coloneqq \left( g \circ f_X\circ Q\right)(t)$ by using the Fa\`a di Bruno chain rule given in \eqref{eq:FaadiBruno}; starting with the latter Taylor expansion,
\begin{align}
  \left(g \circ f_X\circ Q \right)(t+hx) &= \sum_{j = 0}^{s \wedge r - 1} \frac{h^j x^j \left( g \circ f_X\circ Q \right)^{(j)}(t)}{j!}  + \frac{h^{s\wedge r} x^{s \wedge r} \left( g \circ f_X \circ Q \right)^{(s \wedge r)}(t + \tau hx)}{(s \wedge r)!},
\label{eq:compositeTaylorg}
\end{align}
where $\left| \tau \right| < 1$. The intermediate derivatives for $j = 0,1,\ldots, s\wedge r$ are given by \[ \left( g \circ f_X \circ Q \right)^{(j)}(t) = \sum_{\mathbf{k} \in {\cal K}_j } (-1)^{k}k! \left((f_X \circ Q)(t)\right)^{-(k+1)} \prod_{\ell = 1}^j \left( \frac{\left( f_X \circ Q\right)^{(\ell)}(t)}{j!}\right)^{k_\ell} = {\cal O}(1)\]
due to restrictions imposed in \hyperref[AssumptionB]{Assumption (B)}. Similarly, 
\begin{equation}
 \left(\widetilde{f}_\eta \circ Q \right)(t+hx) = \sum_{j = 0}^{s \wedge r - 1} \frac{h^j x^j \left( \widetilde{f}_\eta \circ Q \right)^{(j)}(t)}{j!}  + \frac{h^{s\wedge r} x^{s \wedge r} \left( \widetilde{f}_\eta \circ Q \right)^{(s \wedge r)}(t + \delta hx)}{(s \wedge r)!} \label{eq:compositeITaylor}
\end{equation}
where $\left|\delta\right| \le 1$. Therefore, using \eqref{eq:compositeITaylor} and \eqref{eq:compositeTaylorg} in \eqref{eq:CExpK31} with the vanishing moment condition \eqref{eq:K3moment} implies that,
\begin{align*}
 \CExp{K_3\left( \frac{F(X_i)-t}{h}\right)}{{\cal F}_{i-1}} &= h^{s \wedge r + 1} \Bigg\{ \sum_{j = 0}^{s \wedge r - 1} \sum_{\substack{\ell = 0 \\ j + \ell \ge s \wedge r}}^{s \wedge r - 1} h^{\ell + j - s \wedge r}\frac{p^{(j)}(t) q^{(\ell)}(t) }{j! \ell!} \int_{-1}^1 x^{\ell + j}K_3(x) \, dx \nonumber\\
  & \qquad \quad + \sum_{j = 0 }^{s \wedge r - 1} \frac{p^{(j)}(t)h^j}{(s\wedge r)!j!} \int_{-1}^1 x^{s \wedge r + j} K_3(x) q^{(s \wedge r)} ( t + \tau hx) \,dx \nonumber\\
  & \qquad \quad + \sum_{\ell = 0 }^{s \wedge r - 1} \frac{q^{(\ell)}(t)h^\ell}{(s\wedge r)!\ell!} \int_{-1}^1 x^{s \wedge r + \ell} K_3(x) p^{(s \wedge r)} ( t + \delta hx) \,dx \nonumber\\
  & \qquad \quad + \frac{h^{s \wedge r}}{\left((s\wedge r)!\right)^2} \int_{-1}^1 x^{2 (s \wedge r)} K_3(x) q^{(s \wedge r)}(t + \tau hx)p^{(s \wedge r)} ( t + \delta hx) \,dx\Bigg\}.
\end{align*}
However, by \hyperref[AssumptionB]{Assumption (B)}, $f_\eta^{(j)}$ and $Q$ are Lipschitz continuous for $j = 0,\ldots, s$ and therefore bounded. Consequently $p^{(j)}$ and $q^{(j)}$ are also bounded which means that uniformly in $t$,
\begin{equation}
   \CExp{K_3\left( \frac{F(X_i)-t}{h}\right)}{{\cal F}_{i-1}} < Ch^{s \wedge r + 1} \quad \text{a.s.}
  \label{eq:CExpObound}
\end{equation}
Define, $\widetilde{K_3}(X_{i,i-1}, t) \coloneqq \CExp{K_3\left( \frac{F(X_i)-t}{h}\right)}{{\cal F}_{i-1}}$ and $g(X_{i,i-1},t) \coloneqq \widetilde{K_3}^2(X_{i,i-1}, t) - \Exp \widetilde{K_3}^2(X_{i,i-1}, t)$, then $\Exp g(X_{i,i-1}, t) = 0$ and by Jensen's Inequality $\Exp \widetilde{K_3}(X_{i,i-1}, t)^2 < \infty$. It will be shown by an application of Theorem 1 of \cite{Wu-2007} that $\sum_{i = 1}^n g(X_{i,i-1},t) = {\cal O}_p\left( h^{s \wedge r + 2}n^{1 - \frac{\alpha}{2}}L(n)\right)$. Then, define the physical dependence measure, $\vartheta_i = \sup_{t \in (h,1-h)}\norm{\CExp{g(X_{i,i-1},t)}{{\cal F}_0} - \CExp{g(X_{i,i-1},t)}{{\cal F}_{-1}}}$. To bound $\vartheta_i$, let $\eta_0'$ be an i.i.d. copy of $\eta_0$ and define $X_{i,i-1}^* = X_{i,i-1} - c_i \eta_0 + c_i \eta_0'$ with the associated sigma field ${\cal F}_i^* = \sigma\left( \eta_i, \eta_{i-1}, \ldots, \eta_1, \eta_0', \eta_1, \ldots ; \varepsilon_1, \ldots, \varepsilon_i\right)$. Then by Theorem 1 of \cite{Wu-2005} it was shown that $\vartheta_i \le \sup_{t \in (h,1-h)}\norm{g(X_{i,i-1},t) - g(X_{i,i-1}^*,t)}$. Using this, \eqref{eq:CExpObound} and the Lipschitz property of $f_\eta$ it will be shown that $\vartheta_i < C h^{s\wedge r + 2}i^{-\beta}L(i)$,
\begin{align*}
  \vartheta_i &\le \sup_{t \in (h,1-h)}\norm{g(X_{i,i-1},t) - g(X_{i,i-1}^*,t)} \nonumber\\
  &= \sup_{t \in (h,1-h)}\norm{\left(\widetilde{K_3} \left( X_{i,i-1},t\right) + \widetilde{K_3} \left( X^*_{i,i-1},t\right)\right)\left(\widetilde{K_3} \left( X_{i,i-1},t\right) - \widetilde{K_3} \left( X^*_{i,i-1},t\right)\right)} \nonumber\\
  &\le Ch^{s\wedge r + 1}\sup_{t \in (h,1-h)} \norm{\widetilde{K_3} \left( X_{i,i-1},t\right) - \widetilde{K_3} \left( X^*_{i,i-1},t\right)}\nonumber\\
  &= Ch^{s\wedge r + 1}\sup_{t \in (h,1-h)} \norm{\int_\mathbb{R} K_3 \left( \frac{F(u)- t}{h}\right)\left( f_\eta\left( u - X_{i,i-1}\right) - f_\eta\left( u - X^*_{i,i-1}\right) \right)\, du }\nonumber\\
  &\le Ch^{s\wedge r + 1}\sup_{t \in (h,1-h)} \int_\mathbb{R} \left| K_3 \left( \frac{F(u)- t}{h}\right)\right|\, du \norm{X_{i,i-1}- X^*_{i,i-1}}\nonumber\\
  &\le Ch^{s\wedge r + 2}\norm{\eta_0 - \eta_0'}c_i = C h^{s\wedge r + 2} i^{-\beta}L(i),
\end{align*}
where the last line follows due to the Lipschitz property of $Q$ and the bounded domain of $K_3$. Then by Theorem 1 of \cite{Wu-2007} and Karamata's Theorem, $\norm{\sum_{i = 1}^n g(X_{i,i-1},t)}^2 = {\cal O}\left( h^{2(s\wedge r) + 4}n^{2 -\alpha}L^2(n)\right)$. Using this and \eqref{eq:CExpObound},
\begin{align}
 \frac{1}{nh}\sum_{i = 1}^n \CExp{K_3\left( \frac{F(X_i) - t}{h}\right)}{{\cal F}_{i - 1}}^2 &= \frac{1}{nh}\sum_{i = 1}^n g(X_{i,i-1},t) + \frac{1}{nh} \sum_{i = 1}^n \Exp \widetilde{K}_3^2(X_{i,i-1},t)\nonumber\\
  &= {\cal O}_p\left( h^{2 (s \wedge r ) + 3} n^{-\frac{\alpha}{2}}L(n)\right) + {\cal O}\left( h^{2(s\wedge r)+2}\right)\nonumber\\
  	&=  o_p(1)
  \label{eq:CExp2Opbound}
\end{align}
Then the first term on the RHS of \eqref{eq:condvarB} can be bounded by \eqref{eq:CExp2Opbound} and a similar application of Lemma 1 of \cite{Zhao-Wu-2008},
\begin{align}
  \frac{\mu_F^2(t)}{nh}\sum_{i = 1}^n  &\CExp{\left(K_3\left( \frac{F(X_i)-t}{h}\right)-\CExp{K_3\left( \frac{F(X_i)-t}{h}\right)}{{\cal F_{i-1}}}\right)^2}{{\cal F_{i-1}}}\nonumber\\
  &= \frac{\mu_F^2(t)}{nh}\sum_{i = 1}^n \Bigg\{ \Exp K_3^2\left( \frac{F(X_i)-t}{h}\right) + \CExp{K_3\left( \frac{F(X_i)-t}{h}\right)}{{\cal F}_{i-1}}^2 \nonumber\\
  &\qquad \qquad + 
 \left(\CExp{K_3^2\left( \frac{F(X_i)-t}{h}\right)}{{\cal F_{i-1}}} - \Exp K_3^2\left( \frac{F(X_i)-t}{h}\right) \right)\Bigg\}\nonumber\\
  &= \mu_F^2(t)\int_{- 1}^1 K_3^2(x) \, dx + {\cal O}_p(n^{-\frac{\alpha}{2} }L(n)) + {\cal O}\left(h^2\right)
  \label{eq:condvarB3}
\end{align}
Substituting \eqref{eq:condvarB3} and \eqref{eq:condvarB1} into \eqref{eq:condvarB} implies that,
\begin{equation*}
  \sum_{i=1}^n \CExp{\Delta_i^2(t)}{{\cal F}_{i-1}} \convP 1. 
\end{equation*}
For the Lindeberg condition, let $\epsilon > 0$ and define $A_n = \left\{ |\Delta_1(t)| > \epsilon\right\}$, then similar to the procedure used in the Proof of \autoref{CLTA}, it can be shown that $A_n \to \emptyset$ and the Lindeberg condition holds. Thus by the martingale CLT, \eqref{eq:CLTB} holds and by using \eqref{eq:bandB1} in the decomposition given in \eqref{eq:CLTBdecomp} the result follows by Slutsky's Theorem.
\end{proof}

\begin{proof}[Proof of \autoref{CLTB2}]
Again, use the decomposition \eqref{eq:CLTBdecomp} used in the Proof of \autoref{CLTB1}. Then, define the standardised process,
\[
  \Upsilon_i(t) \coloneqq \frac{\CExp{ \gamma_i(t)}{{\cal F}_{i - 1}}- \Exp \gamma_i(t)}{h^4n^{1- \frac{\alpha}{2}}L(n) {\cal H}_1(t)}.
\]
It will be shown via use of a Hermite expansion of the LRD variables that,
\begin{equation}
	\sum_{i = 1}^n \Upsilon_i(t) \convD {\cal N} (0,C_1^2).
	\label{eq:HermiteCLTB}
\end{equation}
To do this, split the LRD variable $X_i$ into two parts, $X_i = \eta_i + X_{i,i-1}$. Define the standardised version of $X_{i,i-1}$, $Z_i = s_X^{-1}\left(X_{i,i-1}-\mu_X\right)$, $Z_i \sim {\cal N}(0,1)$. Notice that $\Upsilon_i(t)$ and $Z_i$ are both ${\cal F}_{i-1}$-measurable and define $G(Z_i,t) \coloneqq \CExp{\gamma_i(t)}{{\cal F}_{i-1}} - \Exp \gamma_i(t)$. Then clearly, $\Exp G(Z_i,t) = 0$ and by Jensen's inequality, $\Exp G(Z_i,t)^2 < \infty$. So by \cite{Taqqu-1975}, $G(Z_i,t)$ can be re-expressed by its Hermite expansion,\[G(Z_i,t) = \sum_{m = 1}^\infty \frac{a_m}{m!} H_m(Z_i)\] where $a_m = \Exp \left[ H_m(Z_1) G(Z_1,t)\right]$ is the $m^{th}$ Hermite coefficient. For our case it is assumed that $a_1  \neq 0$. Evaluating $a_1$,
\begin{align*}
  a_1&= \Exp\left[Z_1 G(Z_1,t)\right] = \Exp \left[ Z_1\frac{1}{\sigma_\eta}\int_\mathbb{R} \mu(u + \mu_X + s_X Z_1 ) K_3 \left( \frac{\Phi(u + s_X Z_1 ) - t}{h}\right) \phi\left( \frac{u}{\sigma_\eta}\right)\, du\right]\nonumber\\
    &=   \frac{1}{\sigma_\eta}\int_\mathbb{R}\int_\mathbb{R} z\mu(u + \mu_X + s_X z ) K_3 \left( \frac{\Phi(u + s_X z ) - t}{h}\right) \phi(z) \phi\left( \frac{u}{\sigma_\eta}\right)\, dz\,d u\nonumber\\
    &=   \frac{h}{s_X^2 \sigma_\eta}\int_\mathbb{R}\int_{-\frac{t}{h}}^{\frac{1-t}{h}} \frac{\Phi^{-1}(t+hw) - u }{\phi(\Phi^{-1}(t + hw))}\mu_F(t+hw) K_3 \left( w\right) \phi\left(\frac{\Phi^{-1}(t+hw) - u }{s_X}\right) \phi\left( \frac{u}{\sigma_\eta}\right)\, dw\,d u.
 \end{align*}
By exploiting the Fa\`a di Bruno formula further, it can be shown via Taylor expansions that the asymptotic behaviour of $a_1$ satisfies, 
\[ a_1 \sim \frac{h^4 \kappa_h(t)}{s_X^3\sigma_\eta\phi\left( \Phi^{-1} ( t)\right)}\int_\mathbb{R}  \phi\left(\frac{\Phi^{-1}\left( t\right)-u}{s_X}\right) \left(\Phi^{-1}\left( t\right) - u\right)\phi\left( \frac{u}{\sigma_\eta}\right) \, d u= h^4 {\cal H}_1(t) \]
From Corollary 5.1 of \cite{Taqqu-1975},
\begin{align*}
  \sum_{i = 1}^n \Upsilon_i(t) &\sim \frac{1}{n^{1-\frac{\alpha}{2}}L(n)}\sum_{i = 1}^n Z_i \convD {\cal N}(0,C_1^2)
\end{align*}
 Therefore \eqref{eq:HermiteCLTB} holds by Slutsky's Theorem in the decomposition given in \eqref{eq:CLTBdecomp} in conjuction with \eqref{eq:CLTB}, \eqref{eq:HermiteCLTB} and \eqref{eq:bandB2}.
\end{proof}

\begin{proof}[Proof of \autoref{EVTA}]
First, fix $k \in \mathbb{N}$ and choose distinct integers $0 \le j_1, j_2, \ldots, j_k \le m_n$. We adapt the proof of Theorem 5 of \cite{Zhao-Wu-2008} to our case. The proof of their result was reliant on another result given by Theorem 1 of \cite{Grama-Haeusler-2006} which requires a martingale difference sequence that has third order moments. We obtain such a sequence below. Define,
\begin{align*}
	S_{n,k}^A(\mathbf{t}) &=  \sum_{i = 1}^n \Psi_i(\mathbf{t}) = \left[ S_n^A(t_{j_1}), S_n^A(t_{j_2}), \ldots, S_n^A(t_{j_k})\right]^T. 
\end{align*}
and the associated sigma field ${\cal A}_i = \sigma(X_i, \ldots, X_1; \eta_{i+1},\eta_i, \ldots) $. Then $\left\{ S^A_{n,k}(\mathbf{t}), {\cal A}_n\right\}_{i = 1}^n$ is a martingale since $\Exp K_3 \left( \frac{F(X_i)-t}{h}\right) = 0$ for all $t \in (h,1-h)$. Let $\cal{Q}$ be the quadratic characteristic matrix of $S^A_{n,k}$, that is,
\begin{align*}
	\cal{Q} &= \sum_{i = 1}^n \CExp{\Psi_i(\mathbf{t})\Psi_i(\mathbf{t})^T}{{\cal A}_{i-1}} \coloneqq \left( \cal{Q}_{rr'}\right)_{1 \le r, r' \le k}.\\
	\cal{Q}_{rr'}&= \frac{\Exp \left(K_3 \left(\frac{F(X_1) - t_{j_r}}{h}\right)K_3 \left(\frac{F(X_1) - t_{j_{r'}}}{h}\right)\right)	}{n h \,\sqrt{\upsilon(t_{j_{r}})\upsilon(t_{j_{r'}})}}\sum_{i = 1}^n \left(\sigma^2\varepsilon_i^2 + \left( \mu_F(t_{j_r}) + \mu_F(t_{j_{r'}})\right)\sigma\varepsilon_i +  \mu_F(t_{j_r})\mu_F(t_{j_{r'}}) \right).
\end{align*}
However, by construction, if $r \neq r'$, then $|t_{j_r}-t_{j_{r'}}| > 2h$ and the kernel function $K_3 \map{[-1,1]}{\mathbb{R}}$ which implies that $\left\{ x \in \mathbb{R}: \left\{ h^{-1}\left| F(x) - t_{j_r}\right| \le 1 \right\} \cap \left\{ h^{-1}\left| F(x) - t_{j_{r'}}\right| \le 1 \right\}\right\} = \emptyset$. Therefore when $r \neq r'$, $\cal{Q}_{rr'} = 0.$ If $r = r'$, then by \eqref{eq:asympsumvars}, 
\begin{align*}
	\cal{Q}_{rr } &= \frac{1}{n \, \left( \sigma^2 + \mu_F^2(t_{j_r})\right) } \sum_{i = 1}^n \left(\sigma^2\varepsilon_i^2 + 2\mu_F(t_{j_r}) \sigma\varepsilon_i + \mu_F^2(t_{j_r}) \right).\\
	\norm{\cal{Q}_{rr}-1}&\le \frac{1}{n \, \left( \sigma^2 + \mu_F^2(t_{j_r})\right) } \left(\sigma^2\norm{\sum_{i = 1}^n \left(\varepsilon_i^2 -1\right)} + 2 \left| \mu_F(t_{j_r})\right| \norm{\sum_{i = 1}^n \sigma\varepsilon_i} \right)= {\cal O}\left( n^{- \frac{\alpha}{2}} L(n) \right).\nonumber
\end{align*}
Let $\left( u_{rr'} \right)_{1 \le r,r'\le k}$ be the $k \times k$ identity matrix. Then by the above argument $\Exp | \cal{Q}_{rr'}-u_{rr'}|^{\frac{3}{2}} = {\cal O}\left(n^{- \frac{3\alpha}{4}} L^{\frac{3}{2}}(n) \right)$ uniformly over $1 \le r,r' \le k.$ Also, $\sum_{i = 1}^n \Exp |\Psi_i(t)|^3 = {\cal O}(n^{- \frac{1}{2}}h^{- \frac{1}{2}})$. Combining the two yields,  $\sum_{i = 1}^n \Exp |\Psi_i(\mathbf{t})|^3 + \Exp |\cal{Q}_{rr'} - u_{rr'}|^{\frac{3}{2}}= {\cal O}((nh)^{- \frac{1}{2}} + n^{- \frac{3\alpha}{4}}L^{3/2}(n)).$ Considering the asymptotic behaviour of \eqref{eq:Bn}, $\left( 1 + B_{m_n}(x)\right)^4\exp\left\{ \frac{B_{m_n}^2(x)}{2} \right\} = {\cal O}\left( h^{-1} \left|\log h\right|^{\frac{3}{2}} \right)$ and using \eqref{eq:EVTAband} it follows that $\left( 1 + B_{m_n}(x)\right)^4\exp\left\{ \frac{B_{m_n}^2(x)}{2} \right\}\Lambda_n \to 0.$  Therefore the same framework and argument applies that was used in the proof of Theorem 5 of \cite{Zhao-Wu-2008} and the result follows.
\end{proof}

\begin{proof}[Proof of \autoref{EVTB}]
The proof of the Theorem uses a similar result to Theorem 5 of \cite{Zhao-Wu-2008}. However, to be able to adapt the result of \autoref{EVTA} to this case and ensure that $S_n^B(t)$ can be modified into a martingale we add and subtract the conditional expectation by defining,
\begin{equation*}
 S_n^{B^*}(t) \coloneqq \sum_{i = 1}^n \Xi_i^*(t) = \sum_{i = 1}^n \left( \Xi_i(t) - \frac{1}{\upsilon(t)\sqrt{nh}}\CExp{K_3\left( \frac{F(X_i)-t}{h}\right)}{{\cal F}_{i-1}}\right).
\end{equation*}
With this definition, $\left\{ S_n^{B^*}(t),{\cal F}_n \right\}_{n \in \mathbb{Z}^+}$ is a martingale and
\begin{equation}
 S_n^B(t) = S_n^{B^*}(t) + \frac{1}{\upsilon(t)\sqrt{nh}}\sum_{i = 1}^n \CExp{K_3\left( \frac{F(X_i) - t}{h}\right)}{{\cal F}_{i-1}}.
\label{eq:EVTB1}
\end{equation}
The proof of the result will follow from Slutsky's Theorem if the first term on the RHS of \eqref{eq:EVTB1} follows the extreme value distribution and the last term on the RHS of \eqref{eq:EVTB1} converges to zero in probability. From \eqref{eq:CExpObound} and \eqref{eq:EVTBband} it follows that $ n^{-\frac{1}{2}}h^{-\frac{1}{2}}\sum_{i = 1}^n \CExp{K_3\left( \frac{F(X_i) - t}{h}\right)}{{\cal F}_{i-1}} = {\cal O}_{\rm a.s.}\left( n^{\frac{1}{2}}h^{s\wedge r + \frac{1}{2}}\right) = o_{\rm  a.s.}(1)$. Now turn attention to first term on the LHS of \eqref{eq:EVTB1} and apply a similar procedure to the one used in the proof of \autoref{EVTA}. Fix $k \in \mathbb{N}$ and choose distinct integers $0 \le j_1, j_2, \ldots, j_k \le m_n$ and define,
\begin{equation*}
  S^{B^*}_{n,k}(\mathbf{t}) = \left[ S_n^{B^*}(t_{j_1}), S_n^{B^*}(t_{j_2}), \ldots, S_n^{B^*}(t_{j_k})\right]^T. 
\end{equation*}
Then $\left\{ S^{B^*}_{n,k}(\mathbf{t}),{\cal F}_n \right\}_{n \in \mathbb{Z}^+}$ is a martingale. Let $\cal{Q}$ be the quadratic characteristic matrix of $S^{B^*}_{n,k}$, that is,
\begin{align*}
  \cal{Q} &= \sum_{i = 1}^n \CExp{\Xi_i^*(\mathbf{t})\Xi_i^*(\mathbf{t})^T}{{\cal F}_{i-1}} \coloneqq \left( \cal{Q}_{rr'}\right)_{1 \le r, r' \le k}.\nonumber\\
  \cal{Q}_{rr'} &= \sum_{i = 1}^n \CExp{\Xi_i^*(t_{j_r})\Xi_i^*(t_{j_{r'}})}{{\cal F}_{i - 1}}\nonumber\\
  &= \frac{1}{nh\upsilon(t_{j_r}) \upsilon(t_{j_{r'}}) }\sum_{i = 1}^n \Bigg\{ \CExp{\sigma^2(X_i)K_3 \left(\frac{F(X_i) - t_{j_r}}{h}\right)K_3 \left(\frac{F(X_i) - t_{j_{r'}}}{h}\right)}{{\cal F}_{i-1}} \nonumber\\
  &\qquad + \mu_F\left(t_{j_{r}}\right)\mu_F\left(t_{j_{r'}}\right)\CExp{K_3 \left(\frac{F(X_i) - t_{j_r}}{h}\right)K_3 \left(\frac{F(X_i) - t_{j_{r'}}}{h}\right)}{{\cal F}_{i-1}}\nonumber\\
  &\qquad \quad - \mu_F\left(t_{j_{r}}\right)\mu_F\left(t_{j_{r'}}\right)\CExp{K_3 \left(\frac{F(X_i) - t_{j_r}}{h}\right)}{{\cal F}_{i-1}}\CExp{K_3 \left(\frac{F(X_i) - t_{j_{r'}}}{h}\right)}{{\cal F}_{i-1}}\Bigg\}.
\end{align*}
By a similar domain argument that was presented in the proof of \autoref{EVTA}, if $r \neq r'$, $\cal{Q}_{rr'} = 0$. If $r = r'$, then, 
\begin{align}
\cal{Q}_{rr } &= \frac{1}{nh\upsilon^2(t_{j_r})  }\sum_{i = 1}^n\Bigg\{ \CExp{\sigma^2(X_i)K_3^2 \left(\frac{F(X_i) - t_{j_r}}{h}\right)}{{\cal F}_{i-1}} + \mu_F^2\left(t_{j_{r}}\right)\CExp{K_3^2 \left(\frac{F(X_i) - t_{j_r}}{h}\right)}{{\cal F}_{i-1}}\nonumber\\
  &\qquad \qquad \qquad \qquad \qquad \qquad  - \mu_F^2\left(t_{j_{r}}\right)\CExp{K_3 \left(\frac{F(X_i) - t_{j_r}}{h}\right)}{{\cal F}_{i-1}}^2\Bigg\}.
 \label{eq:QrrB}
\end{align}
Therefore, using \eqref{eq:condvarB1} and \eqref{eq:condvarB3} in \eqref{eq:QrrB},
\begin{align}
  \norm{\cal{Q}_{rr} - 1}_{\frac{3}{2}} &\le \norm{\cal{Q}_{rr} - 1} = {\cal O}(\delta)
  \label{eq:EVTB3}
\end{align}
where $\delta = n^{- \frac{\alpha}{2}}L(n) + h^2$. Define $\left(u_{rr'}\right)_{1 \le r, r' \le k}$ to be the $k \times k$ identity matrix , then by \eqref{eq:EVTB3}, uniformly over $r$, $\Exp \left| \cal{Q}_{rr'} - u_{rr'}\right|^{\frac{3}{2}} = {\cal O}(\delta^{\frac{3}{2}}).$ Also, $\sum_{i = 1}^n \Exp \left| \Xi_i^*(t_{j_r})\right|^3 = {\cal O}\left( n^{-\frac{1}{2}} h^{- \frac{1}{2}}\right)$ which implies $\sum_{i = 1}^n \Exp |\Xi_i(t)|^3 + \Exp \left| \cal{Q}_{rr'} - u_{rr'}\right|^{\frac{3}{2}} = {\cal O}(\Lambda_n)$ where $\Lambda_n = n^{-\frac{1}{2}} h^{- \frac{1}{2}} + n^{- \frac{3\alpha}{4}}L^{\frac{3}{2}}(n) + h^3.$
Similarly, due to the bandwidth restriction in \eqref{eq:EVTAband}, $\left( 1 + B_{m_n}(x)\right)^4\exp\left\{ \frac{B_{m_n}^2(x)}{2} \right\}\Lambda_n \to 0$ and by the same argument in the proof of \autoref{EVTA} the result follows.
\end{proof}

\section*{Acknowledgements}
The first author would like to thank his Ph.D. supervisor N. C. Weber for his guidance and support; and Wei Biao Wu for comments which helped in reaching the final result.

\end{document}